\documentclass[lettersize,journal]{IEEEtran}
\usepackage{caption} 
\usepackage{cite}
\usepackage{amsmath,amssymb,amsfonts}
\usepackage[linesnumbered,ruled,vlined]{algorithm2e}
\usepackage{graphicx}
\usepackage{textcomp}
\usepackage{tabularx}
\usepackage{hyperref}
\usepackage{booktabs}
\usepackage{diagbox}
\usepackage{makecell}
\usepackage{hyperref}
\usepackage{stfloats}
\usepackage{subfig}
\usepackage{pifont} % 引入pifont宏包
\usepackage{multirow}
\usepackage{fancyhdr}
\newcommand{\cmark}{$\checkmark$} % 对号
\newcommand{\xmark}{$\times$} % 错号

\newtheorem{remark}{Remark}
\newtheorem{assumption}{Assumption}

\newtheorem{corollary}{Corollary}
\newtheorem{definition}{Definition}
\newtheorem{lemma}{Lemma}
\newtheorem{theorem}{Theorem}
\newtheorem{proof}{Proof}
\renewenvironment{proof}{{\bfseries Proof.}}{}

\fancypagestyle{ieee_notice}{
    \fancyhf{}  % 清除所有页眉页脚
    \fancyhead[C]{\footnotesize This work has been submitted to the IEEE for possible publication.\\
    Copyright may be transferred without notice, after which this version may no longer be accessible.}
}

\begin{document}
\allowdisplaybreaks[4]
\title{Decentralized Nonconvex Robust Optimization over Unsafe Multiagent Systems: System Modeling, Utility, Resilience, and Privacy Analysis}
\author{Jinhui Hu, Guo Chen,~\IEEEmembership{Member,~IEEE}, Huaqing Li,~\IEEEmembership{Senior Member,~IEEE}, Huqiang Cheng,\\ Xiaoyu Guo, and Tingwen Huang,~\IEEEmembership{Fellow,~IEEE}
\thanks{The work was supported in part by the Research Grants Council of Hong Kong under Grant CityU-11210222, in part by the Fundamental Research Funds for the Central Universities under Grant SWU-XDJH202312, in part by the National Natural Science Foundation of China under Grant 62173278, and in part by the Chongqing Science Fund for Distinguished Young Scholars under Grant 2024NSCQ-JQX0103.}
\thanks{
J. Hu is with the School of Automation, Central South University, Changsha 410083, China (e-mail: jinhuihu@csu.edu.cn);
J. Hu and X. Guo are with the Department of Biomedical Engineering, City University of Hong Kong, Kowloon, Hong Kong SAR, China (e-mail: jinhuihu3-c@my.cityu.edu.hk; xiaoyguo@cityu.edu.hk);
G. Chen is with the School of Electrical Engineering and Telecommunications, University of New South Wales, Sydney, NSW 2052, Australia (e-mail: guo.chen@unsw.edu.au);
H. Li is with the Chongqing Key Laboratory of Nonlinear Circuits and Intelligent Information Processing, College of Electronic and Information Engineering, Southwest University, Chongqing 400715, China (e-mail: huaqingli@swu.edu.cn);
H. Cheng is with the Key Laboratory of Dependable Services Computing in Cyber Physical Society-Ministry of Education, College of Computer Science, Chongqing University, Chongqing 400044, China (e-mail: huqiangcheng@126.com);
T. Huang is with the Faculty of Computer Science and Control Engineering, Shenzhen University of Advanced Technology, Shenzhen 518055, China (e-mail: huangtw2024@163.com).
}
}

\markboth{Journal of \LaTeX\ Class Files,~Vol.~XX, No.~XX, XXXX~2024}%
{Shell \MakeLowercase{\textit{et al.}}: A Sample Article Using IEEEtran.cls for IEEE Journals}

% \IEEEpubid{0000--0000/00\$00.00~\copyright~2024 IEEE}

\maketitle
\thispagestyle{ieee_notice}

\begin{abstract}
Privacy leakage and Byzantine issues are two adverse factors to optimization and learning processes of multiagent systems (MASs). Considering an unsafe MAS with these two issues, this paper targets the resolution of a category of nonconvex optimization problems under the Polyak-Łojasiewicz (P-Ł) condition. To address this problem, we first identify and construct the unsafe MAS model. Under this kind of unfavorable MASs, we mask the local gradients with Gaussian noise and adopt a resilient aggregation method, self-centered clipping (SCC), to design a differentially private (DP) and Byzantine-resilient (BR) decentralized stochastic gradient algorithm, dubbed \textit{DP-SCC-PL}, aiming to address a class of nonconvex optimization problems in the presence of both privacy leakage and Byzantine issues. The convergence analysis of \textit{DP-SCC-PL} is challenging, as the convergence error arises from the coupled effects of DP and BR mechanisms, as well as the nonconvex relaxation, which is resolved via seeking the contraction relationships among the disagreement measure of reliable agents before and after the SCC aggregation, together with the optimal gap. Theoretical results not only reveal the trilemma between algorithm utility, resilience, and privacy, but also show that
\textit{DP-SCC-PL} can achieve consensus among all reliable agents. It has also been proven that if there are no privacy issues and Byzantine agents, then the asymptotic exact convergence can be recovered. Numerical experiments verify the utility, resilience, and privacy of \textit{DP-SCC-PL} by tackling a nonconvex optimization problem satisfying the P-Ł condition under various Byzantine attacks.
\end{abstract}


\section{Introduction}
\IEEEPARstart{D}{ecentralized} optimization algorithms (DOAs) play an increasingly pivotal role in the intelligent decision-making process of large-scale MASs \cite{Nedic2018}. Examples for potential applications of DOAs include, but are not limited to, areas such as machine learning \cite{Di2024}, signal processing \cite{Alghunaim2022}, cooperative control \cite{Li2021b}, noncooperative games \cite{Huang2022}, and resource allocation \cite{Ding2024,Yue2025}. The development of MASs is enhanced by DOAs. These algorithms allow agents to perform local updates, distributed computation and storage, and communicate only with their neighbors, thereby preserving individual privacy and mitigating communication bottlenecks inherent in centralized architectures. However, the advancement of MASs also comes with two significant security issues, i.e., users' privacy leakage \cite{Huang2024} and Byzantine agents \cite{Allouah2023}.
\begin{table*}[!h]
\centering
\caption{Comparison of this paper with some related research.}
\renewcommand{\arraystretch}{1.5}
\begin{tabularx}{17.0cm}{cccccc}  % 10cm 減去前兩個欄位寬度後，剩下的通通給
\toprule
\hline                      % 第三欄位使用，文字超出的部份會自動折行
\bf{References}  & \bf{Structure} & \bf{Cost Function} & \bf{Assumption on Boundedness}  & DP & BR  \\
\hline
\cite{Guerraoui2021, Ma2022} & distributed  & strongly-convex and smooth & bounded gradients & \cmark & \cmark\\
\cite{Zhu2022} & distributed & weakly-convex and nonsmooth & bounded gradients & \cmark & \cmark \\
\cite{Karimireddy2021,Allouah2023} & distributed & strongly-convex or nonconvex and smooth & bounded variance and heterogeneity & \cmark & \cmark \\
\cite{Huang2024} &  decentralized & strongly-convex and smooth & not required & \cmark &  \xmark \\
\cite{Wang2023a} &  decentralized & nonconvex and smooth & bounded gradients & \cmark &  \xmark \\
\cite{Hu2023a}  &  decentralized & strongly-convex and composite smooth + nonsmooth & not required & \xmark & \cmark \\
\cite{Ben-ameur2016, Fang2022} &  decentralized & strongly-convex and smooth & not required & \xmark & \cmark \\
\cite{He2022, Wu2023} &  decentralized & nonconvex and smooth & bounded variance and heterogeneity & \xmark & \cmark \\
\cite{Ye2024a} &  decentralized & strongly-convex and nonsmooth & bounded gradients & \cmark & \cmark \\
{\bf{This paper}} & decentralized & nonconvex with the P-Ł condition and smooth  & bounded variance and heterogeneity & \cmark & \cmark \\
\hline
\multicolumn{6}{c}{\makecell[c]{Distributed structure refers to a parameter-server structure.  A smooth function implies that it has Lipschitz continuous gradients.}}\\
\bottomrule
\end{tabularx}
\label{Table 0}
\end{table*}
\subsection{Literature Review}\label{LR}
Differential privacy is a popular strategy to protect users' sensitive information from being disclosed, enabling a mathematical analysis of the privacy guarantees of protected objectives. There are many notable works to achieve differential privacy within a decentralized architecture. To name a few, Huang et al. in \cite{Huang2020} proposed a DP ADMM-type decentralized algorithm via adding Gaussian noise to the decision variable for a class of convex optimization problems. When the objective function is nonconvex, Wang et al. in \cite{Wang2023a} enabled differential privacy for a standard decentralized gradient descent method via masking the local gradients with Gaussian noise. Huang et al. in \cite{Huang2024} proposed a DP decentralized gradient-tracking method through masking local decision variables and gradients with Laplace noise. Through lessening the interaction among agents after convergence, Wang et al. in \cite{Wang2024} designed a decentralized primal-dual algorithm, enabling both differential privacy and exact convergence for a category of convex constrained optimization problems. Wang et al. in \cite{Wang2024a} devised a DP time-varying controller for a multiagent average consensus task via injecting a multiplicative truncated Gaussian noise with time-invariant variance into the state of each agent. However, addressing the privacy issue alone is insufficient, as Byzantine issues pose significant challenges to the consensus and stability of MASs \cite{Wu2023, Gong2024, Koushkbaghi2024}.

Therefore, it is imperative to incorporate resilient aggregation mechanisms into DOAs to mitigate the performance degradation incurred by Byzantine agents. For example, Ben-ameur et al. in \cite{Ben-ameur2016} introduced an idea of norm-penalized approximation based on total variation to achieve Byzantine resilience. Despite that the selection of the penalty parameter poses a challenge to all reliable agents, a superiority of the method \cite{Ben-ameur2016} lies in its relaxed constraints on the network topology. Fang et al. in \cite{Fang2022} designed a screening-based DOA framework, which covers four types of screening mechanisms: coordinate-wise trimmed mean (\textit{CTM}), coordinate-wise median, Krum function, and a combination of Krum and \textit{CTM}. The theoretical analysis is only applicable to the case of \textit{CTM}. He et al. in \cite{He2022} proposed a resilient aggregation mechanism \textit{SCC} via extending \textit{CC} \cite{Karimireddy2021} to a decentralized version for a class of general nonconvex optimization problems, where only first-order stationary points can be attained. Wu et al. in \cite{Wu2023} developed a novel resilient aggregation mechanism \textit{IOS} based on iterative clipping. Recently, Gong et al. in \cite{Gong2024b} proposed a BR distributed controller via trusted-edge mean sequential reduction for a cruise control task.

To date, privacy leakage or Byzantine agents can each be effectively tackled independently in decentralized settings \cite{Huang2020,Wang2023a,Huang2024,Wang2024, Ben-ameur2016,Fang2022,He2022,Wu2023,Koushkbaghi2024,Gong2024b}. The simultaneous presence of these two security issues has received little attention in the decentralized domain, despite the fact that its significance has been recognized by many notable distributed optimization algorithms \cite{Guerraoui2021, Ma2022, Allouah2023, Zhu2022} with a central/master agent for distributed learning tasks. Recent literature  \cite{Ye2024a} designed a DP and BR decentralized algorithmic framework for a class of strongly-convex optimization problems under a bounded-gradient assumption. The theoretical result in \cite{Ye2024a} is inspiring, providing a unified analysis on three resilient aggregation methods under strongly-convex and bounded-gradient assumptions. However, the strongly-convex and bounded-gradient assumptions are stringent and not widely available for many practical problems, such as a least-square problem \cite{Yi2022} and a linear quadratic regulator problem in policy optimization \cite{Fazel2018}, which are nonconvex optimization problems satisfying the P–Ł condition.

\subsection{Motivation and Challenge}\label{MC}
The motivation of this paper is to achieve both differential privacy and Byzantine resilience for decentralized stochastic gradient descent (\textit{DSGD}) based methods, such as those proposed in \cite{Lian2017a, Wang2023a, Wu2023}, without relying on the stringent assumptions of strong convexity and bounded gradients required by recent studies \cite{Ye2024a}. Although differential privacy and Byzantine resilience have been studied individually in recent works \cite{Huang2020,Wang2023a,Huang2024,Fang2022,He2022,Wu2023,Koushkbaghi2024}, the simultaneous analysis of differential privacy and Byzantine resilience within a decentralized nonconvex domain is still non-trivial. This is challenging, as the convergence error arises from the coupled effects of DP and BR mechanisms, as well as the nonconvex relaxation, whose joint impact complicates the analysis.

\subsection{Statement of Contribution}\label{Contri}
The main contributions of this paper are summarized below.
\begin{enumerate}
\item To resolve a class of nonconvex optimization problems under an unsafe MAS in the presence of both privacy leakage and Byzantine issues, this paper designs a DP and BR decentralized stochastic gradient algorithm, dubbed \textit{DP-SCC-PL}. \textit{DP-SCC-PL} simultaneously achieves differential privacy and Byzantine resilience, distinguishing itself from existing decentralized methods that study either DP decentralized methods \cite{Huang2020,Wang2023a,Huang2024,Wang2024,Wang2024a} or BR decentralized methods \cite{Ben-ameur2016,Fang2022,He2022,Wu2023,Koushkbaghi2024}. In contrast to the recent studies \cite{Ye2024a}, \textit{DP-SCC-PL} is not only independent of the stringent bounded-gradient assumption but also proved to be available to a class of nonconvex optimization problems satisfying the P-Ł condition, which finds applications in many practical fields \cite{Fazel2018,Yi2022}.

\item With the aid of the P-Ł condition, this paper proves the optimal gap of the sequence generated by \textit{DP-SCC-PL} in the presence of both privacy leakage and Byzantine issues, while many existing works, such as \cite{He2022,Wu2023}, can only prove the first-order stationary points and fail to characterize the optimal gap. We address the theoretical challenges via seeking the contraction relationships among the disagreement measure of reliable agents before and after aggregation, together with the optimal gap. We have summarized the comparison of this paper with some related research in Table \ref{Table 0}, which highlights the contributions of this paper.

\item As a byproduct, the proposed algorithm achieves guaranteed privacy and utility via injecting Gaussian noise with bounded variance, which can serve as an alternative approach to \cite{Ye2024a} that requires Gaussian noise at a same decaying speed as the employed diminishing step-size.
\end{enumerate}

\subsection{Organization}\label{Organ}
Some preliminaries, including the basic notation, network model and adversary definition, problem formulation, and problem reformulation, are given in Section \ref{Section 2}. Section \ref{Section 3} presents the development and detailed updates of \textit{DP-SCC-PL}. The utility, resilience, and privacy of \textit{DP-SCC-PL} are analyzed in Section \ref{Section 4}. Section \ref{Section 5} performs numerical experiments on a decentralized nonconvex optimization problem satisfying the P-Ł condition to verify the utility, resilience, and privacy of \textit{DP-SCC-PL} under various Byzantine attacks. We draw a conclusion and state our future direction in Section \ref{Section 6}.

\section{Preliminaries}\label{Section 2}
\subsection{Basic Notation}\label{Section 2-1}
All vectors in this paper are assumed to be column vectors unless otherwise specified. The remaining basic symbols and notations of this paper are summarized in Table \ref{Table 1}.
\begin{table}[!h]
\centering
\caption{Basic notations.}
\begin{tabularx}{8.8cm}{lX}  % 10cm 減去前兩個欄位寬度後，剩下的通通給
\toprule                % 第三欄位使用，文字超出的部份會自動折行
\hline
\bf{Symbols}  & \bf{Definitions}  \\
\midrule
${\mathbb R}$, ${{\mathbb R}^n}$, ${{\mathbb R}^{m \times n}}$ & the sets of real numbers, $n$-dimensional column real vectors, $m \times n$ real matrices, respectively\\
:= & the definition symbol\\
${\cdot ^ \top }$ & the transpose of any matrices or vectors \\
%${\rm{diag}}\left\{ \nu \right\}$ & a diagonal matrix with all the elements of vector $\nu \in {\mathbb{R}^n}$ laying on its main diagonal\\
%$X \le Y$ & each element in $Y - X$  is nonnegative, where $X$ and $Y$ are two vectors or matrices with same dimensions\\
% $\tilde x \otimes \tilde y$ & the Kronecker product of vectors $\tilde x$ and $\tilde y$\\
$\left| \cdot \right|$ & an operator to represent the absolute value of a constant or the cardinality of a set\\
$\left\|  \cdot   \right\|_2$ & the standard Euclidean norm for vectors or spectral norm for matrices\\
$\left\|  \cdot   \right\|_F$ & the Frobenius norm for matrices\\
${\mathbf{I}}$   & an identity matrix with appropriate dimension\\
% ${0_{n}}$   & an $n$-dimensional column vector with all-zero elements\\
${\mathbf{1}}$   & an all-ones vector with appropriate dimension\\
$x \sim N\left( {{\tilde \mu }, {{\tilde \sigma }^2}{\mathbf{I}}} \right)$ & to indicate the variable $x$ subject to a Gaussian distribution with expectation $\tilde \mu$ and variance ${{\tilde \sigma }^2}{\mathbf{I}}$ in an element-wise manner\\
\hline
\bottomrule
\end{tabularx}
\label{Table 1}
\end{table}

\subsection{System Model and Adversary Definition}\label{Section 2-2}
\begin{figure}[!h]
  \centering
  \includegraphics[width=3.0in]{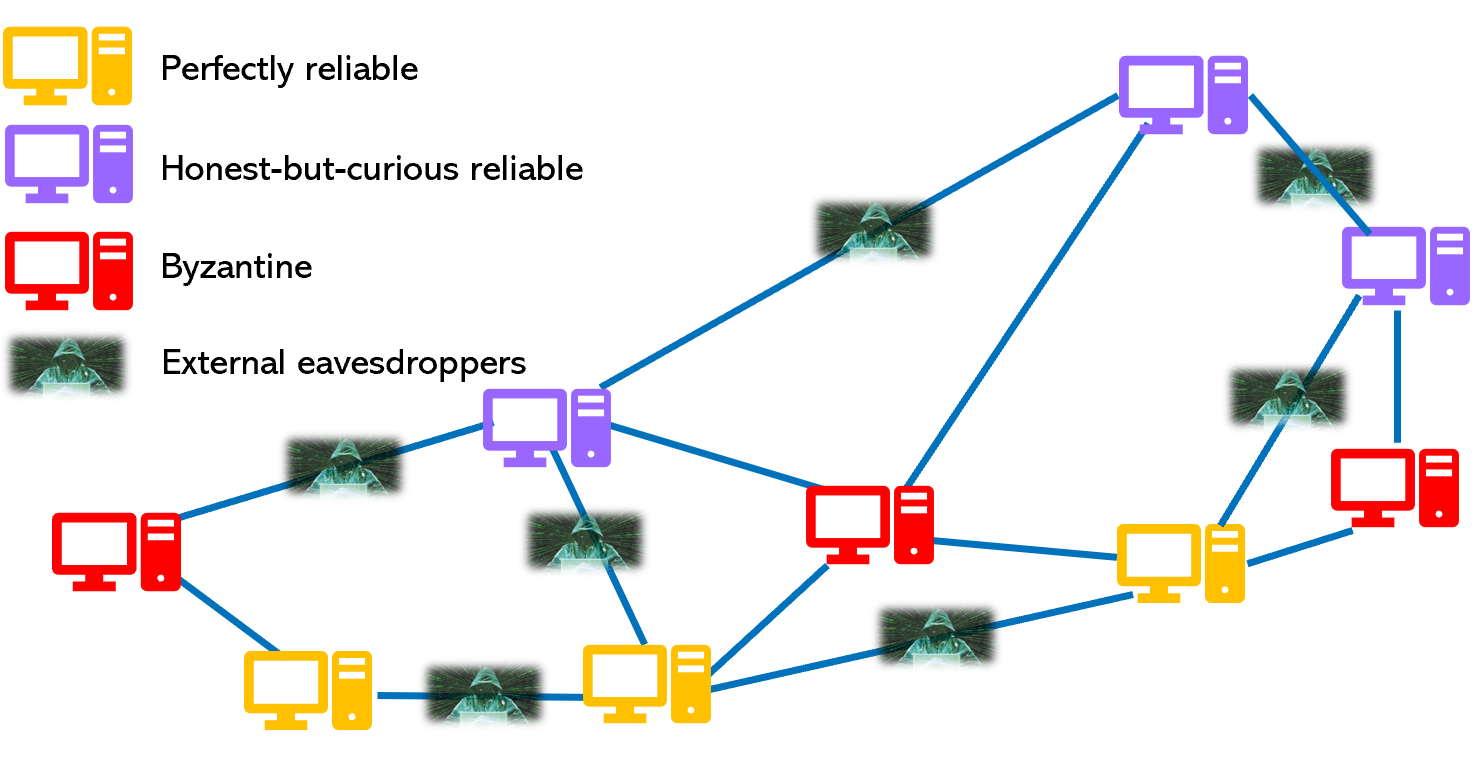}
  \caption{A structural example of an unsafe networked system.}\label{Fig.0}
\end{figure}
We consider a static undirected network $\mathcal{G}: = \left( {\mathcal{V},\mathcal{E}} \right)$ in the presence two kinds of security issues, where $\mathcal{V}$ and $\mathcal{E}$ denotes the set of all agents and communication links over networks, respectively. The first security threat is the existence of Byzantine agents over networks. The sets of reliable and Byzantine agents are denoted by $\mathcal{R}$ and $\mathcal{B}$, respectively. The second threat is the privacy leakage, incurred by two types of adversaries: \textit{honest-but-curious} adversaries and external \textit{eavesdroppers}. Fig. \ref{Fig.0} is an example to briefly describe a MAS consisting of perfectly reliable agents, \textit{honest-but-curious} reliable agents, Byzantine agents, and external \textit{eavesdroppers}. The specific descriptions of Byzantine agents and privacy adversaries are given as follows:
\begin{enumerate}
  \item Byzantine agents are either malfunctioning or malicious agents caused by many possible factors in the course of optimization, such as poisoning data, software bugs, damaged devices, and cyber attacks \cite{Gong2024b}. To study the worst case of the Byzantine model, all Byzantine agents are assumed to be omniscient and able to disobey the prescribed update rules. So, they may collude with each other and send maliciously-falsified information to their reliable neighbors at each iteration \cite{Wang2022}. The impact of Byzantine agents on their reliable neighbors and even the whole MAS has been analyzed by \cite{Hu2023a,Wu2023}.
  \item \textit{Honest-but-curious} adversaries are reliable agents that hold curiosity about some sensitive messages. Therefore, they follow all the update rules to collect all received models and learn the sensitive information about other participants, possibly in a collusive manner. An \textit{honest-but-curious} agent $i$, ${i \in \mathcal{R}}$, has the knowledge of internal information, for instance gradient updates and aggregation rules, but fails to know any messages that are not destined to it \cite{Wang2023a,Wang2024}. Note that \textit{honest-but-curious} agents are not considered Byzantine, as Byzantine agents are assumed to have full knowledge of the entire network.
  \item External \textit{eavesdroppers} are outside adversaries that eavesdrop on communication channels to intercept intermediate messages transferring among agents to learn the sensitive information. So, they have the knowledge of any shared information but fail to get access to any interval information \cite{Wang2023a,Wang2024}. Note that external \textit{eavesdroppers} are different from Byzantine agents, since the latter are internal participants.
\end{enumerate}
 Note that perfectly reliable agents work normally and will not actively introduce any privacy issues. The simultaneous presence of privacy issues and Byzantine agents brings great challenges to the optimization and learning processes process of MASs, since these two issues may not only separately impose a negative influence on the utility \cite{Fang2022,Huang2024} of optimization algorithms but collectively introduce coupling errors \cite{Guerraoui2021, Ma2022, Allouah2023, Zhu2022} to the convergence results.

\begin{assumption}(Network and weight conditions)\label{A1}\\
i) The weight matrix $W: = \left[ {{w_{ij}}} \right]$ associated with $\mathcal{G}$ is nonnegative, i.e., ${w_{ij}} \ge 0$ for $1 \le i,j \le m$, and doubly-stochastic, i.e., $W{\mathbf{1}} = {\mathbf{1}}$ and ${{\mathbf{1}}^ \top }W = {{\mathbf{1}}^ \top }$. In addition, the diagonal weights $w_{ii}$ associated with the reliable agent $i$, $\forall i \in \mathcal{R}$, are positive; \\
ii) All reliable agents form a connected undirected network $ \mathcal{G}_\mathcal{R}: = \left( {\mathcal{R},{\mathcal{E}_\mathcal{R}}} \right) $.
\end{assumption}
\begin{remark}\label{R1-1}
Assumption \ref{A1}-i) keeps in line with the primitive weight condition presumed by decentralized non-resilient optimization algorithms \cite{Alghunaim2022, Xin2022} that require all diagonal weights to be positive since all participants are assumed to be reliable.
Assumption \ref{A1}-ii) is standard in BR decentralized optimization \cite{Ben-ameur2016,Ma2022,Yemini2022a,He2022}, which ensures an information flow between any two reliable agents.
\end{remark}

\subsection{Problem Formulation}\label{Section 2-3}
Considering a MAS suffers from the privacy and Byzantine issues as stated in Section \ref{Section 2-2}, where two unknown sets of reliable and Byzantine agents are denoted as $\mathcal{R}$ and $\mathcal{B}$, respectively. The identities of \textit{honest-but-curious} adversaries and external \textit{eavesdroppers} are also assumed to be unknown and cannot be purged as well. In this adverse scenario, all reliable agents cooperatively minimize
\begin{equation}\label{E2-1-1}
\textbf{P1}: \quad \mathop {\min }\limits_{\tilde x \in {\mathbb{R}^n} } f\left( \tilde x \right): = \frac{1}{{\left| \mathcal{R} \right|}}\sum\limits_{i \in \mathcal{R}} {{f_i}\left( \tilde x \right)},
\end{equation}
where $ \tilde x $ is the decision variable; ${f_i}\left( \tilde x \right): = {\mathbb{E}_{{\xi _i} \sim {\mathcal{D}_i}}}{f_i}\left( {\tilde x,{\xi _i}} \right)$ denotes the local objective function, where ${\xi _i} $ is a random variable subject to a local distribution ${\mathcal{D}_i}$. With a slight abuse of notation, the subsequent analysis briefly uses $\mathbb{E} \cdot $ to denote the expectation of all related variables. To further specify the problem formulation, we need the following assumptions.

\begin{assumption}(Feasibility)\label{A2}
The global objective function has a lower bound ${f^*}: = {\inf _{\tilde x \in {\mathbb{R}^n}}}f\left( x \right)$ such that $- \infty < {f^*} \le f\left( \tilde x \right)$.
\end{assumption}

\begin{assumption}($L_i$-smooth)\label{A3}
Each local objective function $f_i$, ${ i \in \mathcal{R}}$, has Lipschitz gradients such that for any two vectors $\tilde x,\tilde y \in {\mathbb{R}^n}$, it holds
\begin{equation}\label{E2-1-2}
f_i\left( {\tilde x} \right) - f_i\left( {\tilde y} \right) - \left\langle {\nabla f_i\left( {\tilde x} \right),\tilde y - \tilde x} \right\rangle  \le \frac{{{L_i}}}{2}\left\| {\tilde y - \tilde x} \right\|_2^2,
\end{equation}
where $L_i>0$.
\end{assumption}

For simplicity, we denote $L: = {\max _{i \in \mathcal{R}}}{L_i}$ subsequently.

\begin{assumption}(Independent sampling)\label{A4}
The stochastic sampling processes associated with random vector sequences ${\left\{ {{\xi _{i,k}}} \right\}_{i \in \mathcal{R},k \ge 0}}$ are independent of iterations $k = 0.1,2, \ldots $ and all reliable agents. Moreover, the stochastic sampling is uniform and without replacement.
\end{assumption}

\begin{assumption}(Bounded variance and heterogeneity)\label{A5}
For each reliable agent $i$, ${i \in \mathcal{R}}$ and $\forall \tilde x \in {\mathbb{R}^n}$, we have \newline
i) the variance of its stochastic gradients is bounded and there exists a positive constant $ {\sigma }$ such that
\begin{equation}\label{E2-1-3}
{\sigma ^2}: = \mathbb{E}\left\| {\nabla {f_i}\left( {\tilde x,{\xi _i}} \right) - \nabla {f_i}\left( x \right)} \right\|_2^2 < \infty;
\end{equation}
ii) the heterogeneity of its gradients calculated from the distribution ${{\xi _i} \sim {\mathcal{D}_i}}$ is bounded and there exists a positive constant ${\zeta }$ such that
\begin{equation}\label{E2-1-4}
{\zeta ^2}: = \mathop {\max }\limits_{i \in \mathcal{R}} \left\| {\mathbb{E}\nabla {f_i}\left( {\tilde x,{\xi _i}} \right) - \frac{1}{{\left| \mathcal{R} \right|}}\sum\limits_{j \in \mathcal{R}} {\mathbb{E}\nabla {f_j}\left( {\tilde x,{\xi _j}} \right)} } \right\|_2^2 < \infty.
\end{equation}
\end{assumption}
\begin{remark}\label{R1-1}
Assumptions \ref{A2}-\ref{A5} are standard in decentralized stochastic nonconvex optimization \cite{Lian2017a,Liu2020,He2022,Wu2023}. Under Assumption \ref{A3}, it can be verified that the global objective function $f$ is also $L$-smooth. The bounded-gradient assumption assumed by \cite{Wang2023a,Allouah2023} can be a sufficient but not necessary condition for Assumption \ref{A5} in some cases.
\end{remark}
\begin{assumption}(P-Ł condition)\label{A6}
The global objective function $f\left( \tilde x \right)$ satisfies the P-Ł condition such that for a positive constant $\nu $, it holds
\begin{equation}\label{E2-1-5}
\frac{1}{2}\left\| {\nabla f\left( \tilde x \right)} \right\|_2^2 \ge \nu \left( {f\left( \tilde x \right) - {f^*}} \right).
\end{equation}
\end{assumption}
\begin{remark}\label{R1-2}
The P-Ł condition is well-studied by recent literature, such as \cite{Yi2022,Xu2024}. However, these works are confined to an ideal situation that both privacy leakage and Byzantine agents are absent. The development of a decentralized method to counteract these two issues under the P-Ł condition is challenging since its convergence error can be contributed jointly by DP and BR mechanisms, as well as the nonconvex relaxation, which deserves further investigation. The lack of jointly robust mechanisms to counteract these two issues renders decentralized methods vulnerable in practice \cite{Ben-ameur2016, Fang2022, He2022, Wu2023, Huang2020,Wang2023a,Huang2024,Wang2024}.
\end{remark}

\subsection{Problem Reformulation}\label{Section 2-4}
To resolve \textbf{P1} in a decentralized manner, we introduce a matrix $X = {\left[ {{x_1},{x_2}, \ldots ,{x_{\left| \mathcal{R} \right|}}} \right]^ \top } \in {\mathbb{R}^{\left| \mathcal{R} \right| \times n}}$ that collects local copies $x_i$ of the decision variable $\tilde x$ such that \textbf{P1} can be equivalently written into the following formulation
\begin{equation}\label{E2-2-1}
\begin{aligned}
\textbf{P2}: \quad &\mathop {\min }\limits_{X \in {\mathbb{R}^{\left| \mathcal{R} \right| \times n}}} F\left( X  \right): = \frac{1}{{\left| \mathcal{R} \right|}}\sum\limits_{i \in \mathcal{R}} {{f_i}\left( {{x_i}} \right)},\\
&{\text{subject to  }}{x_i} = {x_j}, \left( {i,j} \right) \in {\mathcal{E}_\mathcal{R}},
\end{aligned}
\end{equation}
where ${x_i} = {x_j}$, $\left( {i,j} \right) \in {\mathcal{E}_\mathcal{R}}$, $ i \in \mathcal{R}$, is the consensus constraint.

\section{Algorithm Development}\label{Section 3}
To simultaneously achieve differential privacy and Byzantine resilience for \textit{DSGD}-based methods, such as \cite{Lian2017a,Wu2023,Wang2023a,Ye2024a}, while independent of two stringent assumptions (strong convexity and bounded gradients), we study a resilient aggregation rule \textit{SCC} \cite{He2022}, which is a decentralized version of the centered clipping method \cite{Karimireddy2021}. Compared with \cite{He2022}, we further inject Gaussian noise into local stochastic gradients at each iteration, which guarantees the differential privacy of \textit{DP-SCC-PL} (see Section \ref{Section 4-4}). Different from \cite{He2022,Wang2023a}, both decaying and constant step-sizes are considered in \textit{DP-SCC-PL}, which allows users to make an appropriate choice according to their customized needs. Corresponding comprehensive results regarding these two step-size strategies are provided in Section \ref{Section 4-3}. We next explain the detailed update of \textit{DP-SCC-PL}. For every reliable agent $i$, \textit{SCC} takes its own model denoted by $\tilde x_i^i$, as a self-centered reference to clip the received models denoted by $\tilde x_j^i$, $j \in {\mathcal{N}_i}: = {\mathcal{R}_i} \cup {\mathcal{B}_i}$.
At each iteration, the update rule of \textit{SCC} takes the form of
\begin{equation}\label{E3-1}
\begin{aligned}
SCC_i\left\{ {\tilde x_{i}^i,{{\left\{ {\tilde x_{j}^i} \right\}}_{j \in {\mathcal{N}_i}}}} \right\} =  \sum\limits_{j \in {\mathcal{N}_i}} {{w_{ij}}\left( {\tilde x_{i}^i + Clip\left\{ {\tilde x_{j}^i - \tilde x_{i}^i,{\tau _{i}}} \right\}} \right)},
\end{aligned}
\end{equation}
where $\tau_i$ is a clipping parameter, $Clip\{ {\tilde x_{j}^i - \tilde x_{i}^i,{\tau _{i}}} \}: = ( {\tilde x_{j}^i - \tilde x_{i}^i} ) \cdot \min \{ {1,{\tau _{i}}/{{\| {\tilde x_{j}^i - \tilde x_{i}^i} \|}_2}} \}$, and $w_{ij}$ is the weight assigned by the reliable agent $i$ to its incoming information of the neighboring agent $j$. The detailed updates of \textit{DP-SCC-PL} are presented in Algorithm \ref{Algorithm 1}. Note that the symbol * means an arbitrary vector in ${\mathbb{R}^n}$. If the Byzantine agent $j$ sends nothing at any iteration, then its neighbors $i \in {\mathcal{N}_j}$ set $\tilde x_{j,k}^i = {\mathbf{0}}$ after the synchronous waiting time, where ${\mathbf{0}}$ is an all-zeros vector with appropriate dimension.
\begin{algorithm}[!h]
\small
    \SetKwBlock{Initialize}{Initialize:}{}
    \SetKwBlock{Repeat}{Repeat for $ k = 0,1,2, \dots$}{}
    \SetKwBlock{End}{End for a required criterion is satisfied}{}
    %\newcommand{\End}{\textbf{End}}
    %\SetKwInOut{Output}{Output}
    \KwIn{a proper decaying or constant step-size $\alpha_k$, and additive Gaussian noise ${{\tilde n}_{i,k}} \sim ~ N\left( {{\text{0}},\tilde \varpi_{i,k} ^2{\mathbf{I}}} \right) \in {\mathbb{R}^n}$ with bounded variance $0 < \tilde \varpi _{i,k}^2 \le \varpi _i^2$.}
    \Initialize{Decision variables $x_{i,0} \in {\mathbb{R}^n}$, $i \in \mathcal{R}$.}
    \For{$k=0,1,\ldots,K-1$}{
        \For{\rm{each reliable agent} $i$, $i \in \mathcal{R}$}{
            {\bf{Calculate}} a local stochastic gradient $\nabla {f_i}\left( {{x_{i,k}};{\xi _{i,k}}} \right)$;\\
            {\bf{Mask}} the local gradient with Gaussian noise ${{\tilde n}_{i,k}}$
            \begin{equation*}\label{E3-3}
            {{\tilde g}_i}\left( {{x_{i,k}}} \right) = \nabla {f_i}\left( {{x_{i,k}};{\xi _{i,k}}} \right) + {{\tilde n}_{i,k}};
            \end{equation*}
            {\bf{Execute}} a local gradient descent step
            \begin{equation}\label{E3-3}
            \tilde x_{i,k}^i = {x_{i,k}} - {\alpha _k}{{\tilde g}_i}\left( {{x_{i,k}}} \right);
            \end{equation}\\
            {\bf{Transmit}} $\tilde x_{i,k}^j = \tilde x_{i,k}^i$ to all neighbors $j$, $j \in {\mathcal{N}_i}$;\\
            {\bf{Receive}} $\tilde x_{j,k}^i$ from all neighbors $j$, $j \in {\mathcal{N}_i}$;\\
            {\bf{Aggregate}} the received information according to
            \begin{equation*}\label{E3-4}
            {x_{i,k + 1}} = {SCC_i}\left\{ {\tilde x_{i,k}^i,{{\left\{ {\tilde x_{j,k}^i} \right\}}_{j \in {{\mathcal{N}_i} = {\mathcal{R}_i} \cup {\mathcal{B}_i}}}}} \right\}.
            \end{equation*}
            }
        \For{\rm{each Byzantine agent} $i \in \mathcal{B}$}{{\bf{Send}} $\tilde x_{i,k}^j =  *$ to all neighbors $j$, $j \in {\mathcal{N}_i}$.}
        }
    \KwOut{all decision variables ${x_{i,K}}$, $i \in \mathcal{V}$}
    \caption{\textit{DP-SCC-PL}.}
    \label{Algorithm 1}
\end{algorithm}
% \footnotetext{The symbol * means an arbitrary vector in ${\mathbb{R}^n}$. If the Byzantine agent $j$ sends nothing at any iteration, then its neighbors $i \in {\mathcal{N}_j}$ set $\tilde x_{j,k}^i = {\mathbf{0}}$ after the synchronous waiting time, where ${\mathbf{0}}$ is an all-zeros vector with appropriate dimension.}
\begin{remark}\label{R3-1}
Note that Algorithm \ref{Algorithm 1} performs an information-based rather than identity-based clipping at each iteration, which implies that no agent will be eliminated from the MAS during the optimization and learning process. Even though Algorithm \ref{Algorithm 1} outputs the decision variables of all participants, including both reliable and Byzantine agents, the presence of faulty models or incorrect decision variables produced by Byzantine agents has no impact on the inference or decision-making processes of reliable agents in decentralized MASs.
\end{remark}

\section{Theoretical Analysis}\label{Section 4}
To facilitate the following analysis, we define vectors ${{\bar x}_k}: = \left( {1/\left| \mathcal{R} \right|} \right)\sum\nolimits_{i \in \mathcal{R}} {{x_{i,k}}} $ and ${{\overset{\lower0.5em\hbox{$\smash{\scriptscriptstyle\frown}$}}{x} }_k}: = \left( {1/\left| \mathcal{R} \right|} \right)\sum\nolimits_{i \in \mathcal{R}} {\tilde x_{i,k}^i} $, matrices ${X_k}: = \left[ {{x_{1,k}},{x_{2,k}}, \ldots ,{x_{\left| \mathcal{R} \right|,k}}} \right]^ \top \in {\mathbb{R}^{\left| \mathcal{R} \right| \times n}}$, ${{\tilde X}_k}: = \left[ {\tilde x_{1,k}^1,\tilde x_{2,k}^2, \ldots ,\tilde x_{\left| \mathcal{R} \right|,k}^{\left| \mathcal{R} \right|}} \right]^ \top \in {\mathbb{R}^{\left| \mathcal{R} \right| \times n}}$, and $\nabla F\left( {{X_k}} \right): = {\left[ {\nabla {f_1}\left( {{x_{1,k}}} \right),\nabla {f_2}\left( {{x_{2,k}}} \right), \ldots ,\nabla {f_{\left| \mathcal{R} \right|}}\left( {{x_{\left| \mathcal{R} \right|,k}}} \right)} \right]^ \top } \in {\mathbb{R}^{\left| \mathcal{R} \right| \times n}}$.

\subsection{Sketch of The Proof}\label{Section 4-1}
Let $\mathbb{E}f_{K + 1}^{{\text{best}}}: = {\min _{k  \in \left\{ {1,2, \ldots ,K + 1} \right\}}}f\left( {{{\bar x}_k }} \right)$ with $K \ge 2$.
To analyze the consensus and convergence of \textit{DP-SCC-PL} to the nonconvex optimization problem (\ref{E2-2-1}), we need to seek contraction relationships among the following three terms:
\begin{enumerate}
  \item the disagreement measure of reliable agents before SCC aggregation: $\mathbb{E}{{\tilde D}_k} := \mathbb{E}\left\| {{{\tilde X}_k} - \frac{1}{{\left| \mathcal{R} \right|}}{{\mathbf{1}}}{\mathbf{1}}^ \top {{\tilde X}_k}} \right\|_F^2$;\\
  \item the disagreement measure of reliable agents after SCC aggregation/consensus error: $\mathbb{E}{D_k} := \mathbb{E}\left\| {{X_k} - \frac{1}{{\left| \mathcal{R} \right|}}{{\mathbf{1}}}{\mathbf{1}}^ \top {X_k}} \right\|_F^2$;\\
%  \item the first-order stationary condition: $\left( {1/\left( {K + 1} \right)} \right)\sum\nolimits_{k = 0}^K {\mathbb{E}\left\| {\nabla f\left( {{{\bar x}_k}} \right)} \right\|_2^2} $ for general nonconvex function $f$ (without Assumption \ref{A6}),
  \item the optimal gap: $\mathbb{E}f_{K + 1}^{\text{best}} - {f^*}$ for any function $f$ satisfying the P-Ł condition.
\end{enumerate}

Note that the technical line of the theoretical analysis is different from that in \cite{Ye2024a}, as this paper does not rely on either the strongly-convex or bounded-gradient assumptions.

\subsection{Consensus Analysis}\label{Section 4-2}
\begin{definition}\label{D0}
For a connected network $\mathcal{G}$, if there is a matrix $\tilde W := \left[ {{{\tilde w}_{ij}}} \right] \in {\mathbb{R}^{\left| \mathcal{R} \right| \times \left| \mathcal{R} \right|}} $, whose $\left( {i,j} \right)$-th entry is given by
\begin{equation}\label{E4-1-1}
	{\tilde w_{ij}} = \left\{ \begin{gathered}
		{w_{ii}} + \sum\limits_{j \in {\mathcal{B}_i}} {{w_{ij}}} ,j = i, \hfill \\
		{w_{ij}},j \ne i, \hfill \\
	\end{gathered}  \right.
\end{equation}
and moreover there exists a constant $\rho \ge 0$ satisfying
\begin{equation}\label{E4-1-2}
		{\left\| {{SCC_i}\left( {\tilde x_i^i,{{\left\{ {\tilde x_j^i} \right\}}_{j \in {\mathcal{R}_i} \cup {\mathcal{B}_i}}}} \right) -{{\hat x}_i} } \right\|_2} \le  \rho \mathop {\max }\limits_{j \in {\mathcal{R}_i} \cup \left\{ i \right\}} {\left\| {\tilde x_j^i - {{\hat x}_i} } \right\|_2},
\end{equation}
where ${{\hat x}_i}: = \left( {1/\left| \mathcal{R} \right|} \right)\sum\nolimits_{i \in \mathcal{R}} {{{\tilde w}_{ij}}\tilde x_j^i} $ and ${{\hat x}_{i,k}}: = \left( {1/\left| \mathcal{R} \right|} \right)\sum\nolimits_{i \in \mathcal{R}} {{{\tilde w}_{ij}}\tilde x_{j,k}^i} $, then $\tilde W$ is called as a \textit{virtual} weight matrix associated with $\mathcal{G}$ and $\rho$ denotes a contraction constant of the resilient aggregation SCC.
\end{definition}

Note that the \textit{virtual} weight matrix $\tilde W$ is not involved in the algorithm updates and is only used for theoretical analysis. We then define a mixing rate (squared) $\lambda : = {\left\| {\tilde W - \left( {1/\left| \mathcal{R} \right|} \right){\mathbf{1}}{{\mathbf{1}}^ \top }} \right\|_2^2}$ to facilitate the theoretical analysis.
\begin{corollary}(A feasible choice for $\left( {{\tau},\rho } \right)$)\label{C1}
Suppose that Assumption \ref{A1} holds. If the clipping parameter is chosen as ${\tau _i} = \sqrt {\left( {1/\sum\nolimits_{j \in {\mathcal{B}_i}} {{w_{ij}}} } \right)\sum\nolimits_{j \in {\mathcal{R}_i}} {{w_{ij}}\left\| {\tilde x_i^i - \tilde x_j^i} \right\|_2^2} } $, then the relation (\ref{E4-1-2}) holds true when the contraction constant satisfies $0 \le \rho  \le 4\mathop {\max }\nolimits_{i \in \mathcal{R}} \sqrt {\sum\nolimits_{r \in {\mathcal{R}_i}} {{w_{ir}}} \sum\nolimits_{b \in {\mathcal{B}_i}} {{w_{ib}}} } $.
	\begin{proof}
	See Appendix \ref{Section 7-1}.
	\end{proof}
\end{corollary}
\begin{remark}\label{R4-1-1}
Corollary \ref{C1} provides a theoretical choice of local clipping parameters ${\tau _i}$, $\forall i \in \mathcal{R}$. However, since both the identity and number of Byzantine agents are not assumed to be prior knowledge, determining the clipping parameter ${\tau _i}$ according to Corollary \ref{C1} is challenging in practice. Therefore, we recommend manually tuning this parameter in practice. Besides, there are many other choices of ${\tau _i}$, for instance, ${\tau _i} = \sum\nolimits_{r \in {\mathcal{R}_i}} {{w_{ir}}\left\| {\tilde x_i^i - \tilde x_r^i} \right\|_2^2} $ despite that it will generate a more conservative upper bound for the contraction constant, i.e., $0 \le \rho  \le 2\mathop {\max }\nolimits_{i \in \mathcal{R}} \sqrt {2\left( {1 + {{\left| {{\mathcal{N}_i}} \right|}^2}} \right)\sum\nolimits_{r \in {\mathcal{R}_i}} {{w_{ir}}} } $. Finding the best choice of the pair $\left( {{\tau},\rho } \right)$ with $\tau : = \left[ {{\tau _1},{\tau _2}, \ldots ,{\tau _{\left| \mathcal{R} \right|}}} \right]$ is beyond the scope of this paper. The feasible selection range of $\rho$ provided in Corollary \ref{C1} implies that for all reliable agents $i \in \mathcal{R}$, if the weights ${w_{ib}} = 0$, $\forall b \in {\mathcal{B}_i}$, then the contraction constant $\rho$ must be zero. However, the absence of Byzantine agents in the network is usually a prerequisite for this scenario to occur.
\end{remark}

%\begin{definition}\label{P1}
%if there exists a non-negative matrix $A: = {\left[ {{\tilde w_{ij}}} \right]}$, $\forall i,j \in \mathcal{R}$, satisfying the following four conditions:
%i) $0 < {\tilde w_{ij}} \le 1$ when $j \in {\mathcal{R}_i} \cup \left\{ i \right\}$;
%ii) ${\tilde w_{ij}} = 0$ when $j \notin {\mathcal{R}_i} \cup \left\{ i \right\}$;
%iii) $\sum\nolimits_{j \in {\mathcal{R}_i} \cup \left\{ i \right\}} {{\tilde w_{ij}}}  = 1$;
%iv) for ${{\hat x}_i}: = \sum\nolimits_{j \in {\mathcal{R}_i}} {{w_{ij}}\tilde x_{j,k}^i} $ and $\forall i \in \mathcal{R}$,
%\begin{equation}\label{E4-1-1}
%{\left\| {{\mathcal{A}_i}\left( {\tilde x_i^i,{{\left\{ {\tilde x_j^i} \right\}}_{j \in {\mathcal{R}_i} \cup {\mathcal{B}_i}}}} \right) -{{\hat x}_i} } \right\|_2} \le  \rho \mathop {\max }\limits_{j \in {\mathcal{R}_i}} {\left\| {\tilde x_j^i - {{\hat x}_i} } \right\|_2},
%\end{equation}
%where $\rho  \ge 0$, then $ \rho $ and $A$ are the contraction constant and \textit{virtual} aggregation associated with the resilient aggregation rule $\mathcal{A}\left( {\cdot, \cdot} \right)$, respectively.
%\end{definition}

Let $\varpi^2 : = {\max _{i \in \mathcal{R}}}{\varpi^2 _i}$. 
% We adopt a standard assumption commonly used in DP literature, for instance \cite{Wang2023a}, to ensure the boundedness of $\varpi^2$ as follows:
%	\begin{assumption}(Bounded variance of noises)\label{A7}
%		For $i \in \mathcal{R}$, the injected Gaussian noise ${{\tilde n}_{i,k}}$ at any iteration $k \ge 0$ has a bounded variance, i.e., $\tilde \varpi _{i,k}^2 \le \varpi _i^2 < \infty$.
%\end{assumption}
The following lemma provides a disagreement measure for all reliable agents before SCC aggregation.
\begin{lemma}(Disagreement measure before SCC aggregation)\label{L2}
Suppose that Assumptions \ref{A1} and \ref{A3}-\ref{A5} hold. We have
\begin{equation}\label{E4-1-3}
\begin{aligned}
\mathbb{E}{{\tilde D}_k} \le & \left( {\frac{1}{{1 - \eta }} + \frac{{12\left| \mathcal{R} \right|{L^2}}}{\eta }\alpha _k^2} \right)\mathbb{E}{D_k} + \frac{{8\left| \mathcal{R} \right|\left( {{\sigma ^2} + {\zeta ^2}} \right)}}{\eta }\alpha _k^2\\
 &+ \frac{{2n\left| \mathcal{R} \right| \varpi ^2 }}{\eta }\alpha _k^2.
\end{aligned}
\end{equation}
	\begin{proof}
	See Appendix \ref{Section 7-2}.
	\end{proof}
\end{lemma}

We define $\varphi : = \lambda  - 4\rho \sqrt {\left| \mathcal{R} \right|}$, $\eta : = \varphi /2$, $ \phi: = \varphi /\left( {4 - \varphi } \right)$, $\vartheta : = 4\left| \mathcal{R} \right|\left( {n{\varpi ^2} + 4\left( {{\sigma ^2} + {\zeta ^2}} \right)} \right)/\phi$,  $\theta : = \phi/\left( {4\sqrt 3 L} \right)$, ${k_0} > 2/\phi$, $\underline{\theta}:= \min \left\{ {\theta ,1/\nu} \right\}$, $\iota  \ge{\left( {1 + 1/{k_0}} \right)^2}$, and $\bar \rho : = \lambda /\left( {4\sqrt {\left| {\mathcal{R}} \right|} } \right)$.
\begin{theorem}(Disagreement measure after SCC aggregation)\label{T1}
Suppose that Assumptions \ref{A1} and \ref{A3}-\ref{A5} hold. If the contraction constant satisfies $0 \le \rho  < \bar \rho$ such that the constants meet $\varphi ,\eta ,\phi \in \left( {0,1} \right)$, and the step-size is decaying and chosen as ${\alpha _k} := \theta /\left( {k + {k_0}} \right)$, then there exists
\begin{equation}\label{E4-1-4}
\mathbb{E}{D_k} \le {\left( {1 - \phi } \right)^k}{D_0} + \frac{{2\iota \vartheta {\theta ^2}}}{\phi }\frac{1}{{{{\left( {k + {k_0}} \right)}^2}}}.
\end{equation}

If the step-size is a constant ${\alpha _k} \equiv \alpha$ and satisfies $0 < \alpha  \le \theta$, then there exists
\begin{equation}\label{E4-1-5}
\mathbb{E}{D_k} \le {\left( {1 - \phi } \right)^k}{D_0} + \frac{\vartheta }{\phi }{\alpha ^2}.
\end{equation}
\begin{proof}
See Appendix \ref{Section 7-3}.
\end{proof}
\end{theorem}
\begin{remark}\label{R4-2-2}
The relation (\ref{E4-1-4}) implies that the consensus of all reliable agents is achieved asymptotically when \textit{DP-SCC-PL} employs the decaying step-size. By contrast, the inequality (\ref{E4-1-5}) establishes a fixed disagreement error of all reliable agents when \textit{DP-SCC-PL} employs the constant step-size. This implies that the proposed algorithm with a proper decaying step-size can achieve a smaller disagreement error.
\end{remark}

\subsection{Convergence Analysis}\label{Section 4-3}
We proceed to derive convergence results for Algorithm \ref{Algorithm 1} with both decaying and constant step-sizes by leveraging the results obtained in Lemma \ref{L2} and Theorem \ref{T1}. For simplicity, we define $L: = {\max _{i \in \mathcal{R}}}{L_i}$.
\begin{theorem}(Decaying step-size)\label{T2}
Suppose that Assumptions \ref{A1}-\ref{A6} holds. If the contraction constant satisfies $0 \le \rho  < \lambda /\left( {4\sqrt {\left| \mathcal{R} \right|} } \right)$ such that the constants meet $\varphi ,\eta ,\phi \in \left( {0,1} \right)$, and the decaying step-size is chosen as ${\alpha _k} := \underline{\theta } /\left( {k + {k_0}} \right)$, then for $K \ge 1$ the convergence sequence of Algorithm \ref{Algorithm 1} is characterized by
\begin{equation}\label{E4-1-6}
\begin{aligned}
& {\mathbb{E}f_{K + 1}^{{\text{best}}} - {f^*}}\\
\le & \frac{{\mathbb{E}f\left( {{{\bar x}_0}} \right) - {f^*}}}{{\underline{\theta } \nu \left( {\ln\left( {K + {k_0}} \right) - \ln \left( {{k_0}} \right)} \right)}} + \frac{{\underline{\theta } L{\sigma ^2}\sum\nolimits_{k = 0}^K {\frac{1}{{{{\left( {k + {k_0}} \right)}^2}}}} }}{{\nu \left( {\ln \left( {K + {k_0}} \right) - \ln \left( {{k_0}} \right)} \right)}} \\
& + \frac{{{L^2}}}{\nu }\left( {\frac{{96\left| \mathcal{R} \right|{\rho ^2}}}{\eta } + \frac{1}{{\left| \mathcal{R} \right|}}} \right) \frac{{\sum\nolimits_{k = 0}^K {\frac{1}{{k + {k_0}}}\mathbb{E}{D_k}} }}{{\ln \left( {K + {k_0}} \right) - \ln \left( {{k_0}} \right)}}\\
& + \frac{{8{\rho ^2}}}{{\nu \left( {1 - \eta } \right){{\underline{\theta } }^2}}}\frac{{\sum\nolimits_{k = 0}^K {\left( {k + {k_0}} \right)\mathbb{E}{D_k}} }}{{\ln \left( {K + {k_0}} \right) - \ln \left( {{k_0}} \right)}} \!+\! \frac{{64\left| \mathcal{R} \right|{\rho ^2}}}{{\nu \eta }}\left( {{\sigma ^2} + {\zeta ^2}} \right)\\
& + \frac{{4n}}{\nu }\left( {1 + \frac{{8\left| \mathcal{R} \right|{\rho ^2}}}{\eta }} \right){\varpi ^2}.
\end{aligned}
\end{equation}

If the contraction constant further satisfies $0 \le \rho  < \bar \rho$, then Algorithm \ref{Algorithm 1} with the decaying step-size generates an upper bound for the optimal gap as follows:
	% in the order of $\mathcal{O}\left( {{\rho ^2}\left( {{\sigma ^2} + {\zeta ^2} + {\varpi ^2}} \right)} \right)$, i.e.,
	\begin{equation}\label{E4-1-7}
		\mathop {{\text{lim}}}\limits_{K \to \infty } \mathbb{E}f_{K + 1}^{{\text{best}}} - {f^*} \le \mathcal{O}\left( {{\rho ^2}\left( {{\varpi ^2}  + {\sigma ^2} + {\zeta ^2}} \right)} + {\varpi ^2} \right).
	\end{equation}
	\begin{proof}
		See Appendix \ref{Section 7-4}.
	\end{proof}
\end{theorem}
\begin{remark}\label{R4-3-1}
When adopting a decaying step-size, Theorem \ref{T2} reveals that Algorithm \ref{A1} converges to a fixed error ball around the optimal value at a rate of $\mathcal{O}\left( {1/\ln K} \right)$ since the first four terms at the right-hand-side (RHS) of (\ref{E4-1-6}) diminishes at the rate of $\mathcal{O}\left( {1/\ln K} \right)$. This convergence rate is comparable to the one established in \cite{Zeng2018} for convex optimization problems. The asymptotic convergence error is characterized by the RHS term of (\ref{E4-1-7}), which consists of the (possibly) untrue aggregation ($\rho ^2$) for Byzantine resilience, the injected Gaussian noise with the bounded variance (${\varpi ^2}$) for differential privacy, the bounded variance (${\sigma ^2}$) for the stochastic gradient estimation, and the bounded heterogeneity (${\zeta ^2}$) among local stochastic gradients.
\end{remark}

The following corollary recovers the asymptotic exact convergence for Algorithm \ref{A1} when there are no privacy issues and Byzantine agents.
\begin{corollary}\label{C4-3-1}
	Under the conditions of Theorem \ref{T2}, if $\varpi = \rho = 0$, then we have $\mathop {{\text{lim}}}\nolimits_{K \to \infty } \mathbb{E}f_{K + 1}^{{\text{best}}} = {f^*}$.
	\begin{proof}
		See Appendix \ref{Section 7-5}.
	\end{proof}
\end{corollary}
\begin{theorem}(Constant step-size)\label{T3}
	Suppose that Assumptions \ref{A1}-\ref{A6} holds. If the contraction constant satisfies $0 \le \rho  < \lambda /\left( {4\sqrt {\left| \mathcal{R} \right|} } \right)$ such that the constants meet $\varphi ,\eta ,\phi \in \left( {0,1} \right)$, and the step-size is a constant ${\alpha _k} \equiv \alpha $ satisfying $0 < \alpha \le \underline{\theta } $, then for $K \ge 0$ the convergence sequence of Algorithm \ref{Algorithm 1} is characterized by
	\begin{equation}\label{E4-1-8}
		\begin{aligned}
			& \mathbb{E}f_{K + 1}^{\text{best}} - {f^*}\\
			\le & \frac{{\mathbb{E}f\left( {{{\bar x}_0}} \right) - {f^*}}}{{\nu \alpha \left( {K + 1} \right)}} + \frac{{\frac{{96\left| \mathcal{R} \right|{L^2}{\rho ^2}}}{\eta } + \frac{{{L^2}}}{{\left| \mathcal{R} \right|}} + \frac{{8{\rho ^2}}}{{1 - \eta }}\frac{1}{{{\alpha ^2}}}}}{{\nu \alpha \left( {K + 1} \right)}}\sum\limits_{k = 0}^K {\mathbb{E}{D_k}} \\
			& + \frac{{L{\sigma ^2}}}{\nu }\alpha + {\frac{{64\left| \mathcal{R} \right|{\rho ^2}}}{{\eta \nu }}\left( {{\sigma ^2} + {\zeta ^2}} \right)} + \frac{{4n}}{{\nu}}\left( {1 + \frac{{8\left| \mathcal{R} \right|{\rho ^2}}}{\eta }} \right){\varpi ^2}.
		\end{aligned}
	\end{equation}

If the contraction constant further satisfies $0 \le \rho  < \bar \rho$, then Algorithm \ref{Algorithm 1} with the constant step-size generates an upper bound for the optimal gap as follows:
\begin{equation}\label{E4-1-9}
\begin{aligned}
			\mathop {{\text{lim}}}\limits_{K \to \infty }\mathbb{E}f_{K + 1}^{\text{best}} - {f^*}  \le & \mathcal{O}\left( {{\rho ^2}\left( {{\varpi ^2} + {\sigma ^2} + {\zeta ^2}} \right)} + {\varpi ^2} \right)\! +\! \alpha \mathcal{O}\left( {{\sigma ^2}} \right)\\
			&+{\alpha ^2}\mathcal{O}\left( {{\rho ^2}\left( {{\varpi ^2} + {\sigma ^2} + {\zeta ^2}} \right)} \right).
\end{aligned}
\end{equation}
\begin{proof}
See Appendix \ref{Section 7-6}.
\end{proof}
\end{theorem}
\begin{remark}\label{R4-3-2}
Since the first two terms on the RHS of (\ref{E4-1-8}) diminishes at a rate of $\mathcal{O}\left( {1/K} \right)$, Theorem \ref{T3} implies that \textit{DP-SCC-PL} converges to a fixed error ball around the optimal value at a sublinear convergence rate of $\mathcal{O}\left( {1/K} \right)$ when adopting a constant step-size, which is faster than the convergence rate $\mathcal{O}\left( {1/\ln K} \right)$ with the decaying step-size.
However, when comparing (\ref{E4-1-7}) and (\ref{E4-1-9}), it also reaches a conclusion that \textit{DP-SCC-PL} with the decaying step-size achieves a smaller asymptotic convergence error than with the constant step-size.
%This convergence speed is comparable to the convergence rate $\mathcal{O}\left( {1/K} \right)$ of \textit{DSGD} \cite{Lan2020} for a class of strongly convex problems and faster than the convergence rate $\mathcal{O}\left( {1/\sqrt K } \right)$ of \textit{D-PSGD} \cite{Lian2017a} for general nonconvex problems.
\end{remark}
\begin{remark}\label{R4-3-3}
From Theorems \ref{T2}-\ref{T3}, we observe that there is no explicit trade-off between the algorithm utility—characterized by both the consensus error and optimal gap/asymptotic convergence error—and the number or proportion of Byzantine agents.  The number or proportion of Byzantine agents affects the utility via the contraction constant $\rho$, which is associated with the
network topologies and has no direct relationship with the number of Byzantine agents. However, the trilemma between utility, resilience, and privacy is still disclosed in Theorems \ref{T2}-\ref{T3}.
\end{remark}

\subsection{Privacy Analysis}\label{Section 4-4}
In this section, we leverage a standard definition of $\left( {\epsilon, \delta } \right)$-differential privacy borrowed from \cite{Dwork2013, Wang2023a}, where $\epsilon$ and $\delta$ represent the privacy/utility trade-off and failure probability, respectively. For any DP mechanism, a smaller $\epsilon$ ensures a higher level of privacy at the expense of a bigger convergence error, while a smaller $\delta$ offers a higher successful probability to achieve differential privacy.
\begin{definition}\label{D1}
Let ${\rm{Range}}\left\{ {h} \right\}$ and ${\rm{Prob}}\left\{ {h} \right\}$ denote the range and probability of the function $h$, respectively. The function $h$ is considered to be $\left( {\epsilon, \delta } \right)$-DP if for all $R \subset {\text{Range}}\left\{ h \right\}$ and
two $\Delta $-adjacent inputs $x$ and $x'$, i.e., ${\left\| {x - x'} \right\|_1} \le \Delta$, the following inequality holds
\begin{equation}\label{E4-1-10}
{\text{Prob}}\left\{ {h\left( {x} \right) \in R} \right\} \le {e^\epsilon}{\text{Prob}}\left\{ {h\left( {x'} \right) \in R} \right\} + \delta.
\end{equation}
\end{definition}
In the sequel, we will show that the injected Gaussian noise can provide both local (single-iteration) and global (end-to-end) differential privacy protection for Algorithm \ref{Algorithm 1}. The global differential privacy relies on an additional requirement of bounded local gradients, which is standard in the relevant literature \cite{Wang2023a}.
Let ${g_{i,k}}: = \nabla {f_i}\left( {{x_{i,k}},{\xi _{i,k}}} \right)$ and ${g'_{i,k}}: = \nabla {f'_i}\left( {{x_{i,k}},{\xi _{i,k}}} \right)$, we define the sensitivity function as follows:
\begin{equation}\label{E4-1-10+}
{S_{i,g}}: = \mathop {\sup }\limits_{{{\left\| {{g_{i,k}} - {g'_{i,k}}} \right\|}_1} \le {{\Delta _g}} } {\left\| {{{\mathcal{A}}_{i,k}}\left( {{g_{i,k}}} \right) - {{\mathcal{A}}_{i,k}}\left( {{g'_{i,k}}} \right)} \right\|_1},
\end{equation}
where $g_{i,k}$ and $g'_{i,k}$ are two $\Delta_g$-adjacent input gradients; ${{{\mathcal{A}}_{i,k}}}$ denotes the update ${{\mathcal{A}}_{i,k}}: = {x_{i,k}} - {\alpha _k}\nabla {f_i}\left( {{x_{i,k}};{\xi _{i,k}}} \right)$. Note that the weights $w_{ij}$, $\forall i,j \in \mathcal{R} \cup \mathcal{B}$, are assumed to be public information, which can be accessed by both \textit{honest-but-curious} adversaries and external \textit{eavesdroppers}.

\begin{theorem}(Local differential privacy)\label{T4}
For any pair of $\left( {\epsilon, \delta } \right)$ with $ 0< \epsilon, \delta < 1$, if each reliable agent $i$, $i \in \mathcal{R}$ employing a decaying step-size ${\alpha _k}$ and the variance ${\varpi ^2}$ satisfies
\begin{equation}\label{E4-1-11}
{\varpi ^2} \ge 2\frac{{{{\Delta_g ^2}{\underline{\theta } }^2}}}{{{k_0^2{\epsilon^2}}}}\left( {\ln \left( {1.25} \right) - \ln \left( \delta  \right)} \right),
\end{equation}
or employing a constant step-size $\alpha$ and the variance ${\varpi ^2}$ satisfies
\begin{equation}\label{E4-1-12}
{\varpi ^2} \ge 2\frac{{{{\Delta_g ^2}{\underline{\theta } }^2}}}{{{\epsilon^2}}}\left( {\ln \left( {1.25} \right) - \ln \left( \delta  \right)} \right),
\end{equation}
then the injected Gaussian noise ${{\tilde n}_{i,k}}$ can ensure $\left( {\epsilon, \delta } \right)$-differential privacy for the local stochastic gradient $\nabla {f_i}\left( {{x_{i,k}};{\xi _{i,k}}} \right)$ at each iteration $k$, $\forall k \ge 0$.
\begin{proof}
See Appendix \ref{Section 7-7}.
\end{proof}
\end{theorem}
\begin{remark}\label{R4-1}
Theorem \ref{T4} implies that a higher level (smaller $\epsilon$)  of DP for local stochastic gradients requires a larger noise variance ${\varpi ^2}$. However, a larger ${\varpi ^2}$ causes a bigger optimal gap/convergence error according to Theorems \ref{T2}-\ref{T3}. Therefore, it demonstrates that there is a trade-off between utility and the level of differential privacy. If a further DP guarantee for training data is required, it is clear from \cite[Theorem 4]{Wang2023a} that need an extra assumption on the Lipschitz continuity of the gradients $\nabla {f_i}\left( {{x};{\xi _{i,k}}} \right)$, $i \in {\cal R}$, for any vector $x$.
\end{remark}
\begin{theorem}(Global differential privacy)\label{T5}
Suppose that all local gradients have a unified upper bound $B_g$, i.e., ${\left\| {\nabla {f_i}\left( {x;{\xi _{i,k}}} \right)} \right\|_2} \le {B_g}$, for any $x \in {{\mathbb R}^n}$,  $\forall i \in {\mathcal{R}}$, $\forall k \ge 0$.
%$\nabla {f_i}\left( {{x};{\xi _{i,k}}} \right)$, $\forall i \in {\mathcal{R}}$, are $\tilde L$-Lipschitz continuous for any vector $x$.
Let $Q$ and $B_s$ denote the number of total samples and the batch size for stochastic gradient sampling, respectively. Given a total number of iterations $K$, if the noise variance ${\varpi ^2}$ and the R\'enyi divergence of order $o_r$ satisfy ${\varpi ^2} \ge 6B_g^2/B_s^2$ and ${o_r} \le  - \log \left( {{B_s}\left( {{\varpi ^2}{Q^2}/\left( {4B_g^2} \right) + 1} \right)/Q} \right)$, respectively, then Algorithm \ref{Algorithm 1} preserves $\left( {\epsilon, \delta } \right)$-differential privacy for all local gradients, where the privacy loss is defined by $\epsilon: = 20B_g^2K/\left( {{\varpi ^2}{Q^2}} \right) + 2{B_g}\sqrt {20K\log \left( {\frac{1}{\delta }} \right)} /\left( {\varpi Q} \right) $ for any probability $\delta > 0$.
\begin{proof}
The proof follows directly from the proof of \cite[Theorem 1]{Xu2022a}.
\end{proof}
\end{theorem}

\begin{figure*}[!h]
%\begin{center}
\subfloat[An star network with $90$ reliable agents and $10$ Byzantine agents.]{\includegraphics[width=1.7in,height=1.2in]{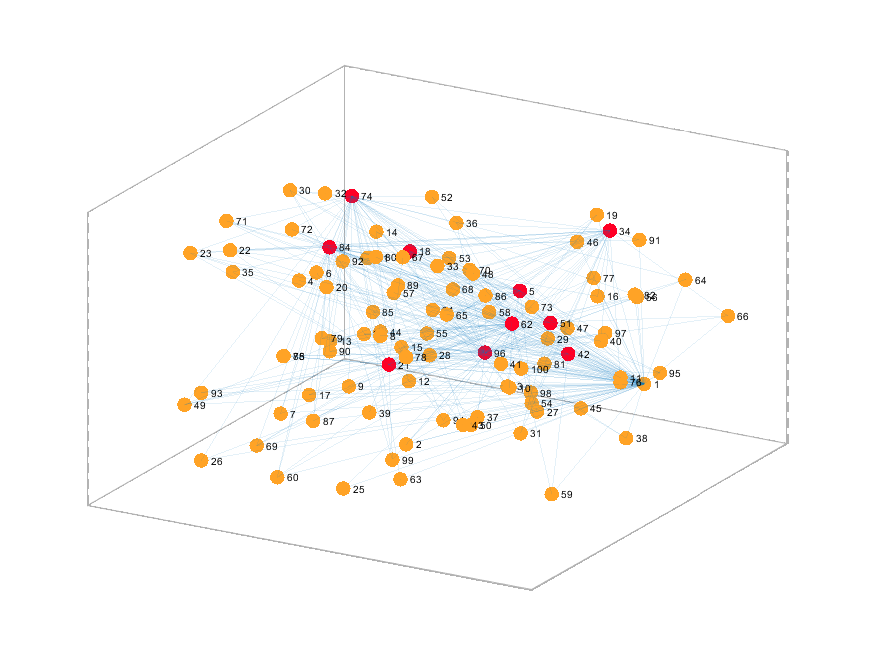}\label{Fig1-1}} \hfill
\subfloat[Consensus error over iterations.]{\includegraphics[width=1.7in,height=1.2in]{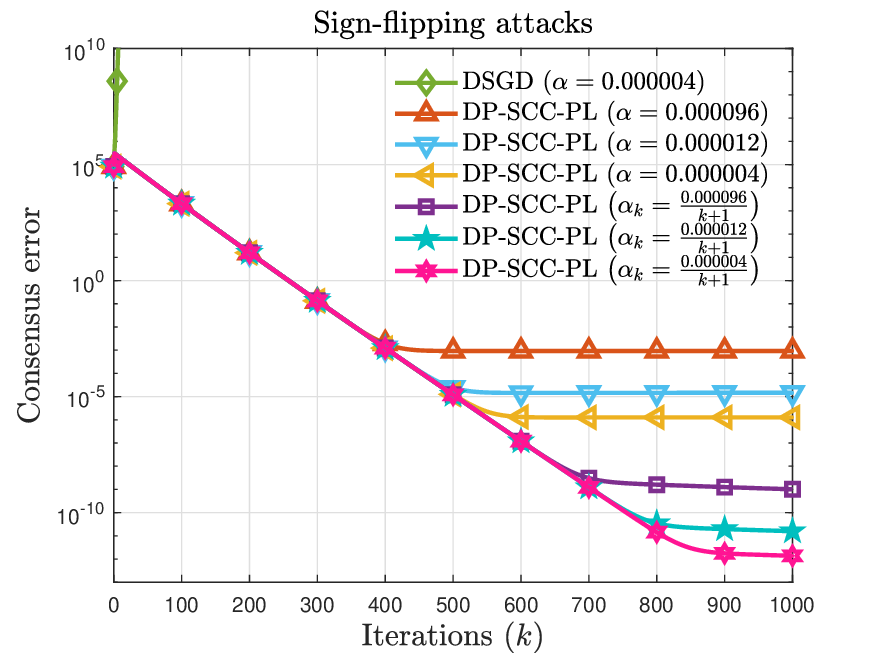}\label{Fig1-2}} \hfill
\subfloat[Optimal gap over iterations.]{\includegraphics[width=1.7in,height=1.2in]{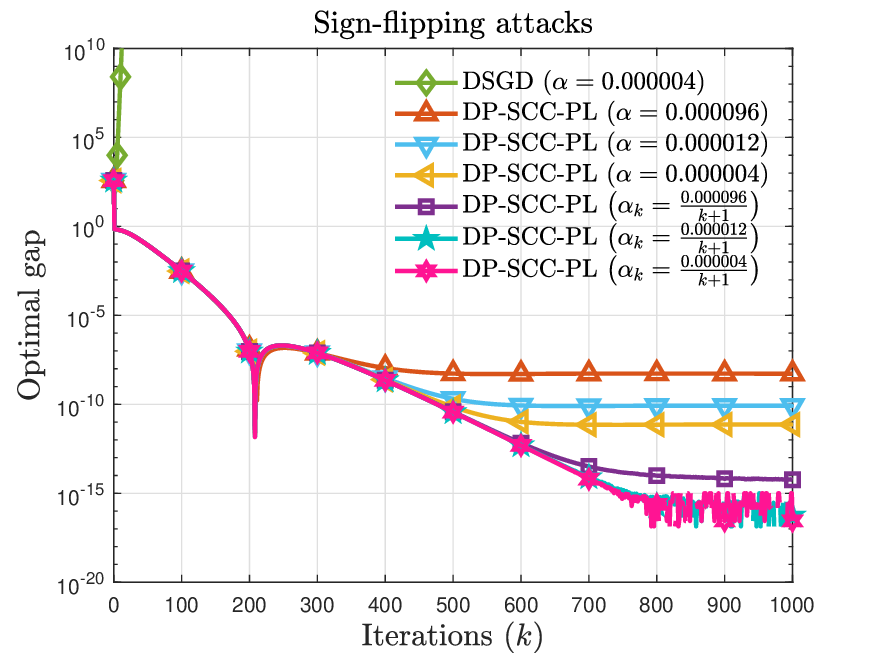}\label{Fig1-3}} \hfill
\subfloat[Observations over iterations.]{\includegraphics[width=1.7in,height=1.2in]{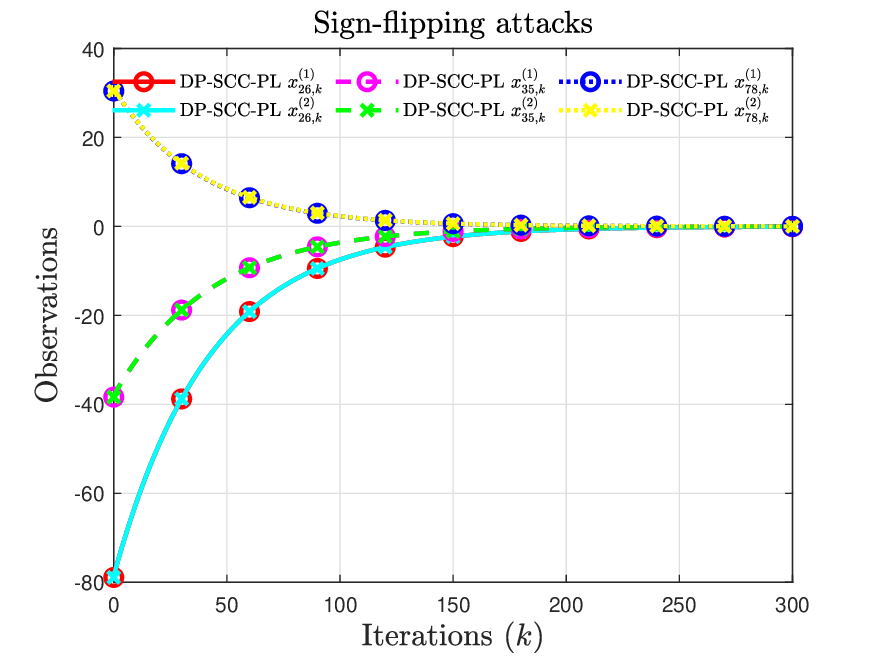}\label{Fig1-4}} \hfill
%\end{center}
\caption{Performance comparison under sign-flipping attacks with injected noises ${{\tilde n}_{i,k}} \sim N\left( {0,0.001^2} \right)$.}
\label{Fig1}
\end{figure*}

\begin{figure*}[!h]
%\begin{center}
\subfloat[An random network with $80$ reliable agents and $20$ Byzantine agents.]{\includegraphics[width=1.7in,height=1.2in]{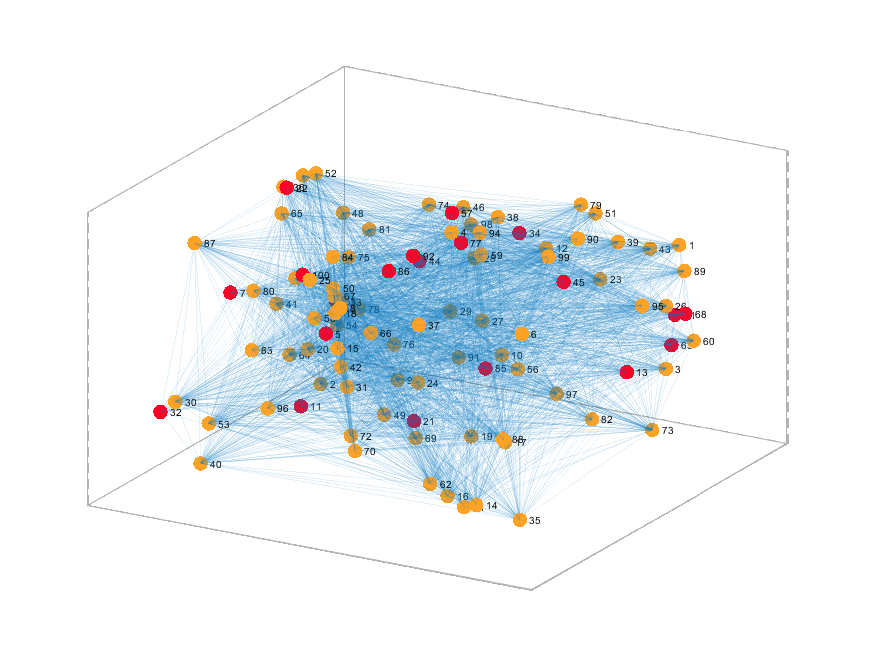}\label{Fig2-1}} \hfill
\subfloat[Consensus error over iterations.]{\includegraphics[width=1.7in,height=1.2in]{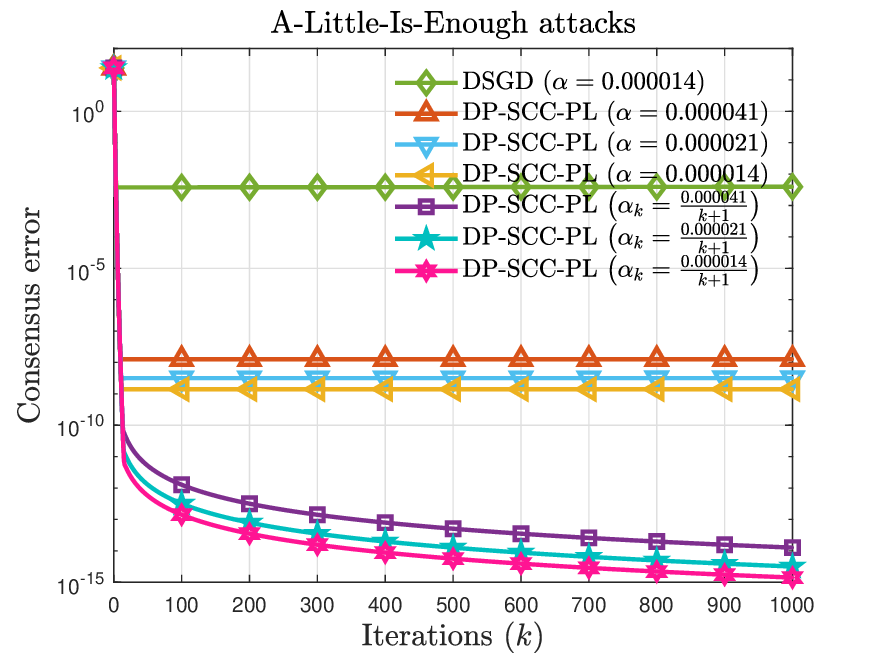}\label{Fig2-2}} \hfill
\subfloat[Optimal gap over iterations.]{\includegraphics[width=1.7in,height=1.2in]{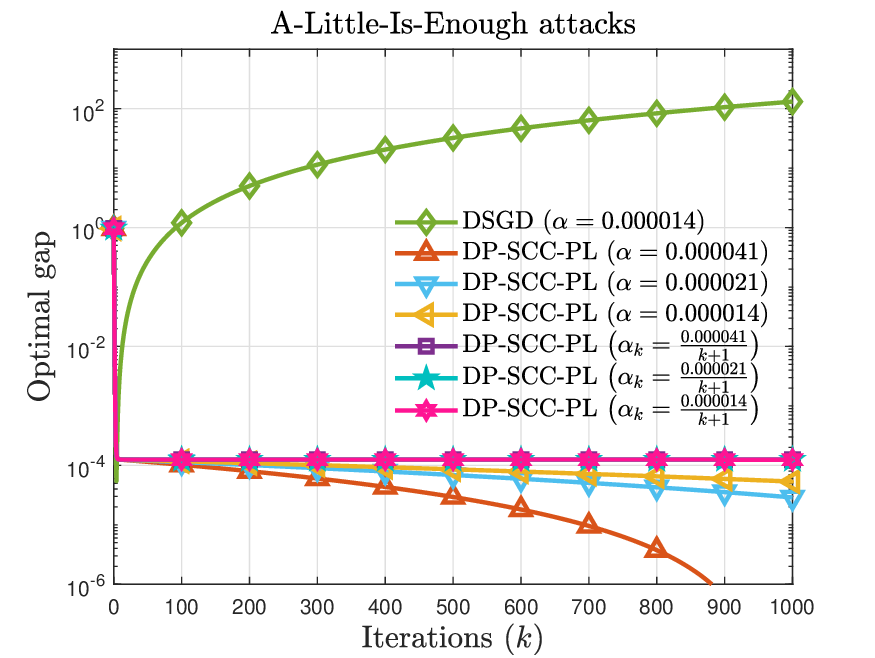}\label{Fig2-3}} \hfill
\subfloat[Observations over iterations.]{\includegraphics[width=1.7in,height=1.2in]{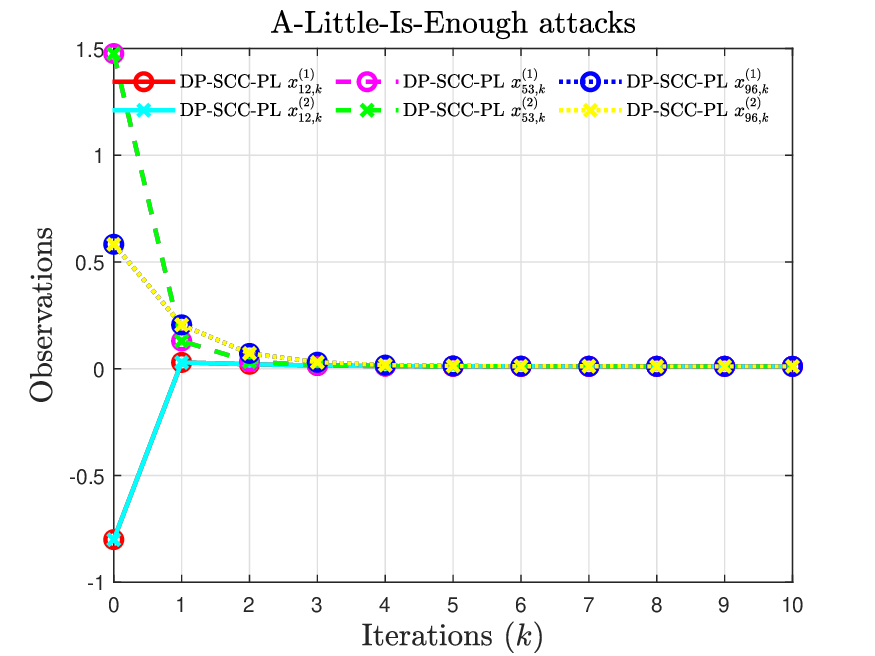}\label{Fig2-4}} \hfill
%\end{center}
\caption{Performance comparison under A-Little-Is-Enough attacks with injected noises ${{\tilde n}_{i,k}} \sim N\left( {0,0.01^2} \right)$.}
\label{Fig2}
\end{figure*}
\begin{figure*}[!h]
%\begin{center}
\subfloat[An fully-connected network with $70$ reliable agents and $30$ Byzantine agents.]{\includegraphics[width=1.7in,height=1.2in]{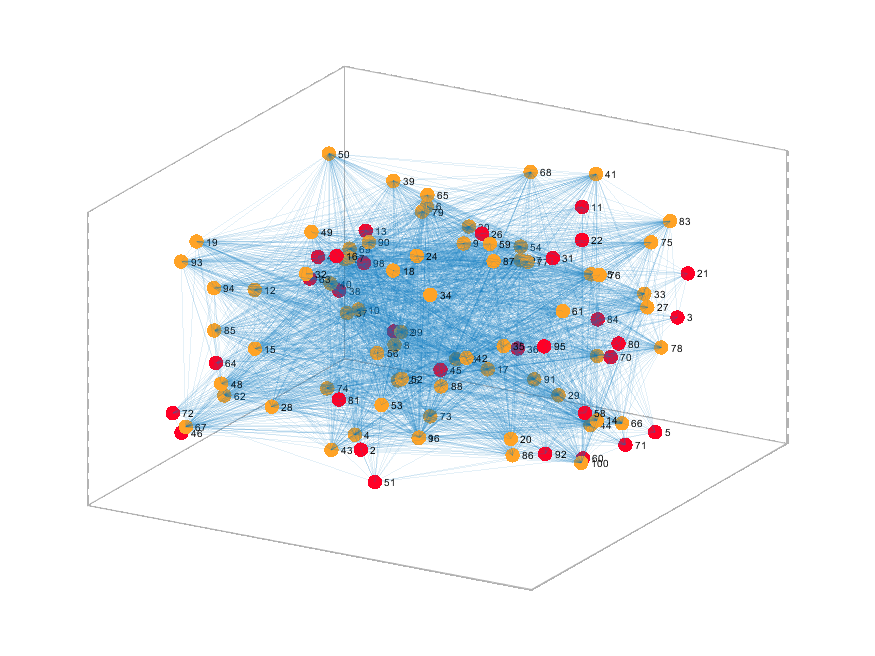}\label{Fig3-1}} \hfill
\subfloat[Consensus error over iterations.]{\includegraphics[width=1.7in,height=1.2in]{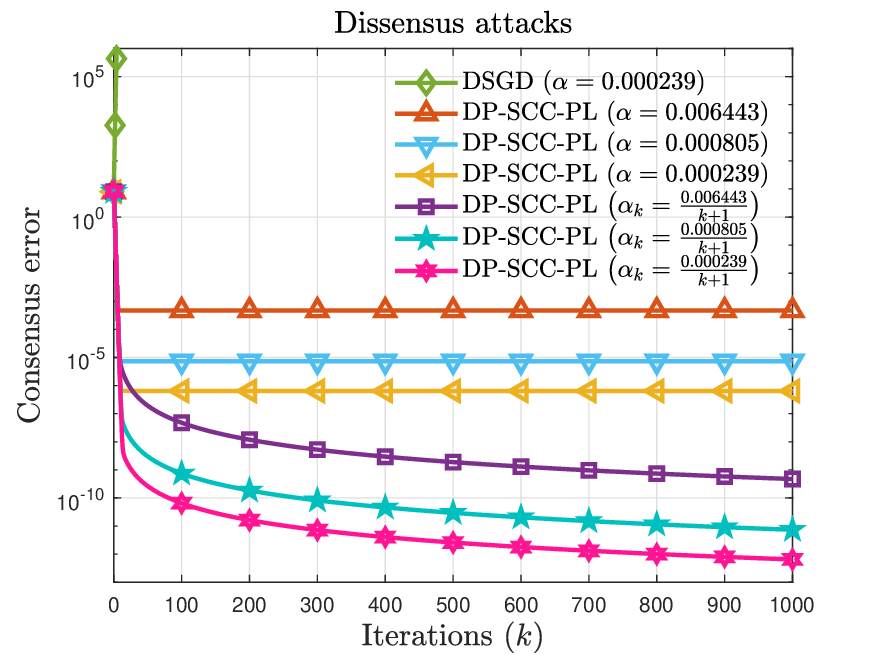}\label{Fig3-2}} \hfill
\subfloat[Optimal gap over iterations.]{\includegraphics[width=1.7in,height=1.2in]{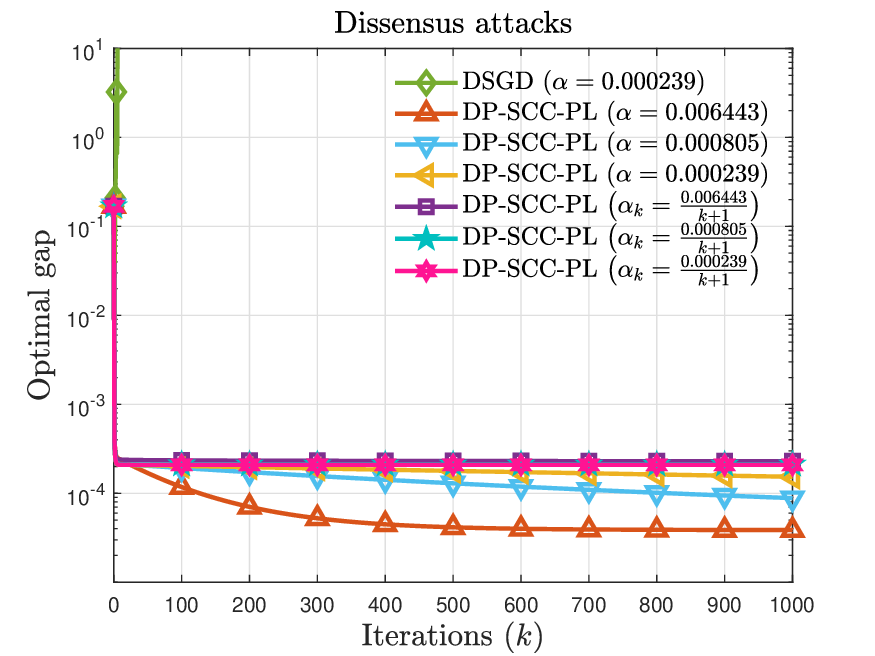}\label{Fig3-3}} \hfill
\subfloat[Observations over iterations.]{\includegraphics[width=1.7in,height=1.2in]{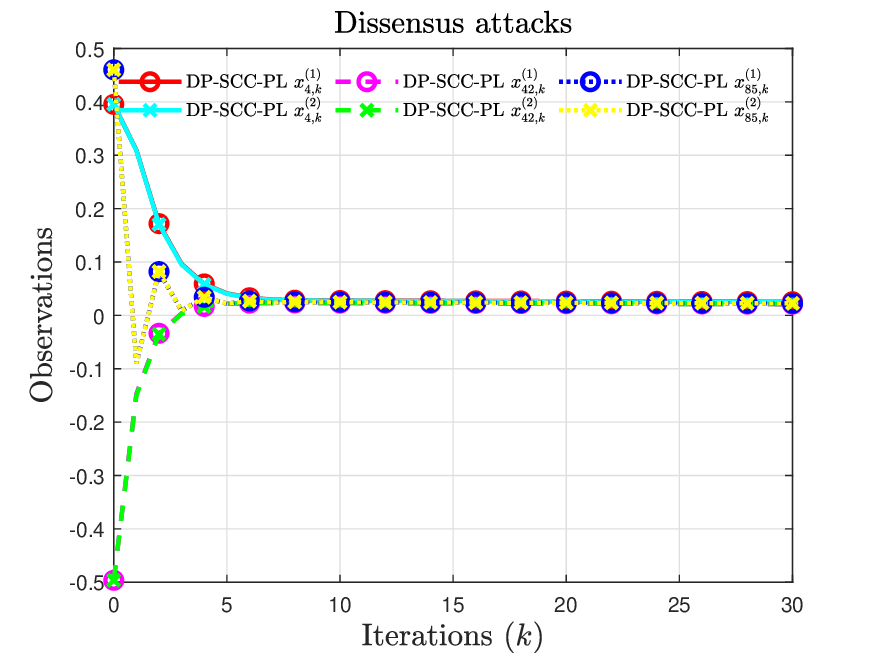}\label{Fig3-4}} \hfill
%\end{center}
\caption{Performance comparison under dissensus attacks with injected noises ${{\tilde n}_{i,k}} \sim N\left( {0,0.1^2} \right)$.}
\label{Fig3}
\end{figure*}

\section{Numerical Experiments}\label{Section 5}
To verify the utility, resilience, and privacy of \textit{DP-SCC-PL}, it is applied to resolving a nonconvex optimization problem over a connected undirected network. A network of $100$ agents are allocated with the sequel local objective functions
\begin{equation*}
\begin{aligned}
{f_j}\left( x \right): = & {{\mathbb{E}_{\left( {{u_j},{v_j}} \right)}}} {0.2u_j\sqrt {{x^4} + 3}  + 0.7u_j{{\cos }^2}x } + u_j,\\
{f_{10 + j}}\left( x \right) := & {\mathbb{E}_{\left( {{u_j},{v_j}} \right)}}{2u_j\sin x - 0.1u_j{{\left( {{x^2} + 2} \right)}^{\frac{1}{3}}} + v_j},\\
{f_{20 + j}}\left( x \right) := & {\mathbb{E}_{\left( {{u_j},{v_j}} \right)}}{\frac{{0.3u_j{x^2}}}{{\sqrt {{x^2} + 1} }} + v_j},\\
{f_{30 + j}}\left( x \right) := & {\mathbb{E}_{\left( {{u_j},{v_j}} \right)}}{ v_j - 0.1u_j\sqrt {{x^4} + 3}  - u_j\sin x},\\
{f_{40 + j}}\left( x \right) := & {\mathbb{E}_{\left( {{u_j},{v_j}} \right)}}{ v_j - \frac{{0.2u_j{x^2}}}{{\sqrt {{x^2} + 1} }} + 2u_j{{\sin }^2}x},\\
{f_{50 + j}}\left( x \right): = & {\mathbb{E}_{\left( {{u_j},{v_j}} \right)}}{ v_j- 0.1u_j\sqrt {{x^4} + 3}  - \frac{{0.1u_j{x^2}}}{{\sqrt {{x^2} + 1} }} },\\
{f_{60 + j}}\left( x \right): = & {\mathbb{E}_{\left( {{u_j},{v_j}} \right)}}{ v_j- u_j\sin x - u_j }, \\
{f_{70 + j}}\left( x \right): = & {\mathbb{E}_{\left( {{u_j},{v_j}} \right)}}{u_j{x^2} + 0.3u_j{{\cos }^2}x + v_j},\\
{f_{80 + j}}\left( x \right): = & {\mathbb{E}_{\left( {{u_j},{v_j}} \right)}}{2u_j{{\sin }^2}x + 0.2u_j{{\left( {{x^2} + 2} \right)}^{\frac{1}{3}}} + v_j},\\
{f_{90 + j}}\left( x \right): = & {\mathbb{E}_{\left( {{u_j},{v_j}} \right)}} {v_j - 0.1u_j{{\left( {{x^2} + 2} \right)}^{\frac{1}{3}}} },
\end{aligned}
\end{equation*}
where $j = 1,2, \ldots ,10$, $ u_j \sim N\left( {1,0.01} \right)$ and $v_j \sim N\left( {0,0.01} \right)$ are two random variables subject to the normal distributions. We denote the function set $\mathcal{F} = {\left\{ {{f_i}} \right\}_{i = 1,2, \ldots ,100}}$. It can be verified that the sum of these local objective functions, i.e., $\sum\nolimits_{i = 1}^{100} {\mathbb{E}{f_i}\left( x \right)}  = 10{x^2} + 30{\sin ^2}x + 10$, is nonconvex but satisfies the P-Ł condition. To ensure the sum of local objective functions of all reliable agents satisfying the P-Ł condition, we evenly choose Byzantine agents from $\left\{ {1,2, \ldots ,100} \right\}$.
\subsection{Case Study One}\label{Section 5-1}
In the first case study, we aim to study the utility, resilience, and privacy of Algorithm \ref{A1} under various scenarios. To assess the utility, we evaluate the consensus error and optimal gap of the proposed algorithm. To verify the differential privacy, the superscripts $\left( 1 \right)$ and $\left( 2 \right)$ are utilized to distinguish the models ${x^{\left( 1 \right)}}$ and ${x^{\left( 2 \right)}}$ with respect to two adjacent function sets ${\mathcal{F}^{\left( 1 \right)}}: = {\left\{ {f_i^{\left( 1 \right)}} \right\}_{i \in \mathcal{R}}} =\mathcal{F} $ and ${\mathcal{F}^{\left( 2 \right)}}: = {\left\{ {f_i^{\left( 2 \right)}} \right\}_{i \in \mathcal{R}}}$, respectively. We randomly choose one function $f_{i_0}$ associated with agent $i_0$ such that it is different between ${\mathcal{F}^{\left( 1 \right)}}$ and ${\mathcal{F}^{\left( 2 \right)}}$ each time while the rest objective functions of ${\mathcal{F}^{\left( 2 \right)}}$ keep same with ${\mathcal{F}^{\left( 1 \right)}}$. To show the Byzantine resilience, we take the following popular Byzantine attacks into consideration.\newline
{\bf{Sign-flipping attacks:}} For any reliable agent $r$, $r \in \mathcal{R}$, its Byzantine neighbors $b \in {\mathcal{B}_r}$ send the falsified model $\tilde x_{b,k}^r =  - {s_b}\sum\nolimits_{i \in {\mathcal{R}_r} \cup \left\{ r \right\}} {{x_{i,k}}} /\left( {\left| {{\mathcal{R}_r}} \right| + 1} \right)$ to it, where $s_b > 0$  is the hyperparameter controlling the attack deviation.\newline
{\bf{A-Little-Is-Enough attacks\cite{Baruch2019}:}} For any reliable agent $r$, $r \in \mathcal{R}$, its Byzantine neighbors $b \in {\mathcal{B}_r}$ send the falsified model $\tilde x_{b,k}^r = {\mu _{{\mathcal{N}_r}}} - a{\sigma _{{\mathcal{N}_r}}}$ to it, where ${\mu _{{\mathcal{N}_r}}}$ and ${\sigma _{{\mathcal{N}_r}}}$ denote the mean and standard deviation of all reliable agents' models, respectively; $a$ is the hyperparameter defined as $a := {\max _a}\left( {\overset{\lower0.5em\hbox{$\smash{\scriptscriptstyle\frown}$}}{c} \left( a \right) < \left( {\left( {\left| \mathcal{V} \right| - \left\lfloor {\left| \mathcal{V} \right|/2 + 1} \right\rfloor } \right)/\left| \mathcal{R} \right|} \right)} \right)$; ${\overset{\lower0.5em\hbox{$\smash{\scriptscriptstyle\frown}$}}{c} }$ is the cumulative standard normal function.\newline
{\bf{Dissensus attacks \cite{He2022}:}} For any reliable agent $r$, $r \in \mathcal{R}$, its Byzantine neighbors $b \in {\mathcal{B}_r}$ send the falsified model $\tilde x_{b,k}^r = {x_{r,k}} - {d_r}\sum\nolimits_{i \in {\mathcal{R}_r}} {{w_{ri}}\left( {{x_{i,k}} - {x_{r,k}}} \right)} /\left( {\sum\nolimits_{b \in {\mathcal{B}_r}} {{w_{rb}}} } \right)$ to it, where $d _r$ is the hyperparameter controlling the attack degree.
% \textit{DSGD} represents the standard gossip-based decentralized stochastic gradient descent algorithm \cite{Lian2017a,Wang2023a}, which is injected with the same Gaussian noise.
Based on the above aforementioned Byzantine attacks, we study the transient behaviors of Algorithm \ref{A1} over three (``star", ``random", and ``full-connected" formulated by reliable agents) classes of undirected networks under different proportions of Byzantine agents and Gaussian noise injection. The decaying and constant step-sizes are selected subject to theoretical bounds ${\alpha _k}: = \underline{\theta } /\left( {k + {k_0}} \right)$ and $\alpha  \in \left( {0,\underline{\theta } } \right]$. Fig. \ref{Fig1} shows that \textit{DP-SCC-PL} with the decaying step-sizes achieves a smaller consensus error and optimal gap than that with the constant step-sizes. In Figs. \ref{Fig2}-\ref{Fig3}, there is a similar outcome that \textit{DP-SCC-PL} with the decaying step-sizes achieves a smaller consensus error than its constant counterpart, whereas \textit{DP-SCC-PL} with the constant step-sizes converges faster than that with decaying step-sizes. From Figs. \ref{Fig1}: \subref{Fig1-4}, \ref{Fig2}: \subref{Fig2-4}, and \ref{Fig3}: \subref{Fig3-4}, we can see that the difference between the models ${\mathcal{F}^{\left( 1 \right)}}$ and ${\mathcal{F}^{\left( 2 \right)}}$ generated from two adjacent function sets in these three case studies is small and almost unobservable. This verifies the differential privacy of \textit{DP-SCC-PL}. Via comparing with the benchmark gossip-based \textit{DSGD} methods \cite{Lian2017a,Wang2023a}, the resilience of \textit{DP-SCC-PL} is verified under various Byzantine attacks. In a nutshell, even though both Gaussian noise and Byzantine attacks are considered, \textit{DP-SCC-PL} still achieves guaranteed consensus and convergence, which is in line with theoretical findings in Theorems \ref{T1}-\ref{T3}.

\subsection{Case Study Two}\label{Section 5-2}
In the second case study, we explore the influence of various proportions of Byzantine agents and different levels of differential privacy on the performance/utility of Algorithm \ref{Algorithm 1}.
\begin{table*}[!htp]
\centering
\caption{Performance comparison across various proportions of Byzantine agents.}
\label{Table 2}
%\scalebox{0.9}{
\begin{tabular}{cccccccc}
\toprule
\hline
\diagbox{Utility}{Proportions} & 0 & 0.1 & 0.2 & 0.3 & 0.4 & 0.5 \\
\midrule
Decaying step-size & 10.8563/($k$+10) & 10.1886/($k$+10) & 50.1338/($k$+10)& 97.4995/($k$+100) & 25.3769/($k$+100) & 23.874/($k$+100)
\\
Consensus error & 4.5249e-11 & {\bf{1.1332e-12}}  &  1.4325e-10  & 2.0940e-09  & 1.1821e-09 & 6.7624e-06\\
Optimal gap & {\bf{7.3571e-08}} & 1.1151e-07 & 1.4498e-07  & 1.8515e-07 &  4.4542e-07 & 8.0112e-04\\
\cline{1-1}
Constant step-size & 5.4281e-03 & 1.0188e-02 & 0.626673 & 0.9749  & 2.5376e-03 & 2.3874e-03 \\
Consensus error & 2.8213e-10  & {\bf{2.8450e-11}} &  1.0670e-06 &  0.9430 & 4.0831e-04 & 6.7840e-04
\\
Optimal gap &  {\bf{7.4027e-08}}  & 1.1172e-07   & 8.2847e-07  & 0.1110 & 1.4159 & 1.4147\\
\hline
%\multicolumn{7}{c}{\makecell[c]{S\&P is an abbreviation of setting and performance.}}\\
\bottomrule
\end{tabular}
%}
\end{table*}
%We consider a network of $m=50$ (including both reliable and Byzantine) agents, which are allocated with local objective functions from part of the function set $\mathcal{F}$ to ensure the sum of them is nonconvex but satisfying P-Ł condition as well.
A novel Byzantine attack is considered in the sequel: \newline
{\bf{Perturbed duplicating attacks:}} For any reliable agent $r$, $r \in \mathcal{R}$, its Byzantine neighbor $b$, $b \in {\mathcal{B}_r}$, sends a same untrue model $\tilde x_{b,k}^r = {\tilde p_{b,k}}{x_{i,k}} + {p_{b,k}}$, $i \in {\mathcal{R}_r}$ to it, where ${\tilde p_{b,k}}$ and $ {p_{b,k}}$ are the multiplicative and additive perturbed parameters, respectively. This kind of attacks can better circumvent the screening, filtration, or clipping processes via selecting appropriate ${\tilde p_{b,k}}$ and $ {p_{b,k}}$. In this experiment, we hand-tune all the parameters, including clipping parameters and step-sizes, to optimize the performance.
%This case study considers a fully-connected network among all agents.
To explore the influence of various proportions of Byzantine agents on the performance of Algorithm \ref{Algorithm 1}, we consider a fixed injection of the Gaussian noise ${{\tilde n}_{i,k}} \sim N\left( {0,0.001^2} \right)$ and five different proportions (0, 0.1, 0.2, 0.3, 0.4, and 0.5) of Byzantine agents in the network.
Although Table \ref{Table 2} shows that as the proportion of Byzantine agents in the network increases, the performance of Algorithm \ref{Algorithm 1} generally deteriorates, evidenced by an increase in both the consensus error and optimal gap, there are exceptions. For example, the consensus errors at the proportion of 0.1 are smaller than that of at the proportion of 0.
%Besides, the smallest optimal gap and consensus error happen at the proportions of 0 and 0.3, respectively. These exceptions exactly verify the statement given in Remark \ref{R4-3-3}, which is in line with theoretical findings in Theorems \ref{T1}-\ref{T3}.
\begin{figure}[!h]
%\begin{center}
\subfloat[Consensus error.]{\includegraphics[width=1.7in,height=1.2in]{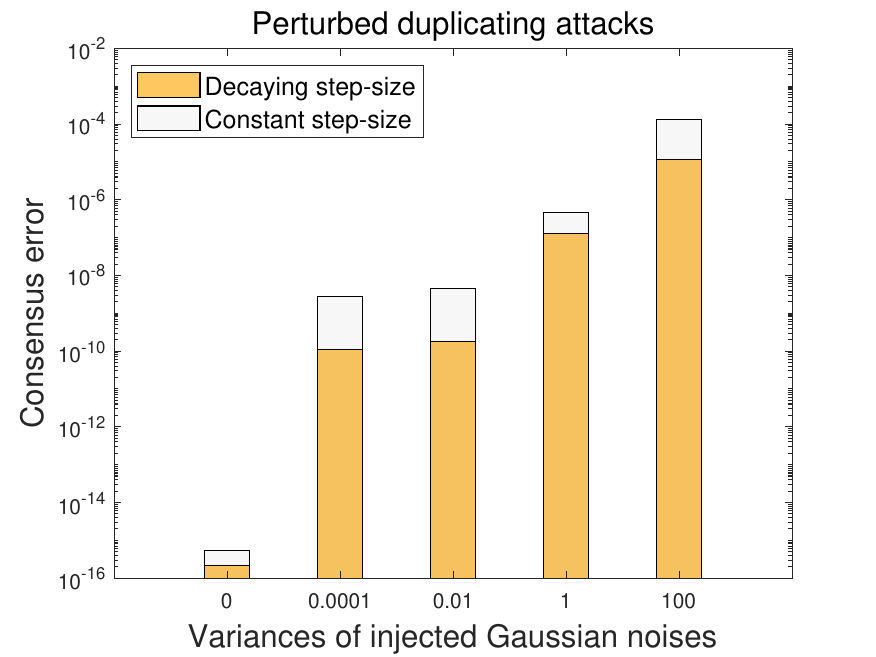}\label{Fig4-1}} \hfill
\subfloat[Optimal gap.]{\includegraphics[width=1.7in,height=1.2in]{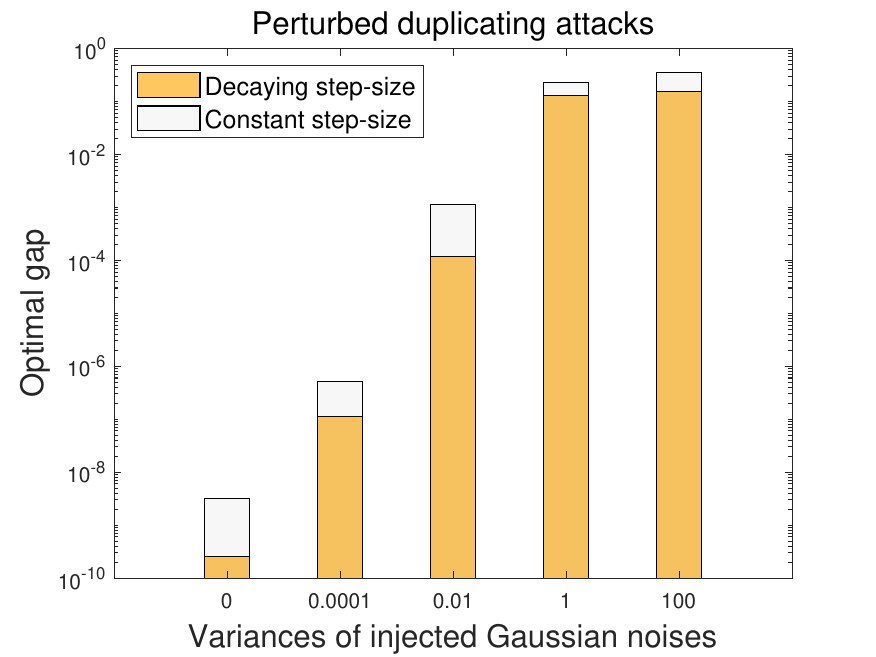}\label{Fig4-2}} \hfill
%\end{center}
\caption{Performance across varying Gaussian noise variance levels.}
\label{Fig4}
\end{figure}
To further explore the influence of different levels of differential privacy on the performance of Algorithm \ref{Algorithm 1}, we consider the injection of Gaussian noise with five kinds of variances ($\varpi^2 = 0, 0.0001, 0.01, 1, 100$). Note that Gaussian noise with larger variance indicates a higher level of differential privacy according to Theorem \ref{T4}. The proportion of Byzantine agents is set as 0.1.
From Figs. \ref{Fig4}: \subref{Fig4-1}-\subref{Fig4-2} show that both the consensus error and the optimality gap exhibit a positive correlation with the magnitude of the noise variance, which verifies the trade-off between utility and differential privacy. From Table \ref{Table 2} and Fig. \ref{Fig4}, it can also be verified that the proposed algorithm with a proper decaying step-size achieves better performance than that with the constant step-size. In a nutshell, these experimental results are consistent with the theoretical findings obtained in Theorems \ref{T1}-\ref{T4}.

\section{Conclusion}\label{Section 6}
This paper studied a nonconvex optimization problem over an unsafe MAS in the presence of both privacy issues and Byzantine agents. To enhance agents' privacy and resilience in the course of optimization, we developed a DP and BR decentralized stochastic gradient algorithm, dubbed \textit{DP-SCC-PL}, via introducing both Gaussian noise and a resilient aggregation mechanism. We addressed the challenge in analyzing the convergence of \textit{DP-SCC-PL} via seeking the contraction relationships among the disagreement measure of reliable agents before and after SCC aggregation, together with the optimal gap. The theoretical results not only reveal the trilemma between utility, resilience, and privacy, but also establish a sublinear (inexact) convergence for \textit{DP-SCC-PL} with a well-designed constant step-size. Numerical experiments verify the utility, resilience, and privacy of \textit{DP-SCC-PL} under various Byzantine attacks via resolving a nonconvex optimization problem satisfying the P-Ł condition. Future directions will extend \textit{DP-SCC-PL} to asynchronous updates.

\section{Appendix}\label{Section 7}
\subsection{Proof of Corollary \ref{C1}}\label{Section 7-1}
For each reliable agent $i$, $i \in \mathcal{R}$, we denote $\tilde z_j^i: = \tilde x_i^i + Clip\left\{ {\tilde x_j^i - \tilde x_i^i,{\tau _i}} \right\}$ and recall the relation (\ref{E3-1}) such that
\begin{equation*}
\begin{aligned}
&\left\| {SC{C_i}\left( {\tilde x_i^i,{{\left\{ {\tilde x_j^i} \right\}}_{j \in {\mathcal{R}_i} \cup {\mathcal{B}_i}}}} \right) - {{\hat x}_i}} \right\|_2^2\\
= & \left\| {\sum\limits_{j \in {\mathcal{N}_i} \cup \left\{ i \right\}} {{w_{ij}}\tilde z_j^i}  - \sum\limits_{j \in {\mathcal{R}_i} \cup \left\{ i \right\}} {{{\tilde w}_{ij}}\tilde x_j^i} } \right\|_2^2\\
= & \left\| {\sum\limits_{j \in {\mathcal{R}_i}} {{w_{ij}}\left( {\tilde z_j^i - \tilde x_j^i} \right)}  + \sum\limits_{j \in {\mathcal{B}_i}} {{w_{ij}}\left( {\tilde z_j^i - \tilde x_i^i} \right)} } \right\|_2^2\\
\end{aligned}
\end{equation*}
\begin{equation}\label{AE7-1-1}
\begin{aligned}
\le & 2\left\| {\sum\limits_{j \in {\mathcal{R}_i}} {{w_{ij}}\left( {\tilde z_j^i - \tilde x_j^i} \right)} } \right\|_2^2 + 2\left\| {\sum\limits_{j \in {\mathcal{B}_i}} {{w_{ij}}\left( {\tilde z_j^i - \tilde x_i^i} \right)} } \right\|_2^2,
\end{aligned}
\end{equation}
where the second equality is according to (\ref{E4-1-1}) and $\tilde z_i^i=\tilde x_i^i$. An upper bound for $\left\| {\sum\nolimits_{j \in {\mathcal{R}_i}} {{w_{ij}}\left( {\tilde z_j^i - \tilde x_j^i} \right)} } \right\|_2^2$ can be derived as follows:
\begin{equation}\label{AE7-1-2}
\left\| {\sum\limits_{j \in {\mathcal{R}_i}} {{w_{ij}}\left( {\tilde z_j^i - \tilde x_j^i} \right)} } \right\|_2^2 \le {\left( {\frac{1}{{{\tau _i}}}\sum\limits_{j \in {\mathcal{R}_i}} {{w_{ij}}{{\left\| {\tilde x_i^i - \tilde x_j^i} \right\|}_2^2}} } \right)^2},
\end{equation}
which applies the fact that $\tilde z_j^i - \tilde x_j^i = 0$ if no clipping happens and ${\left\| {\tilde x_i^i - \tilde x_j^i} \right\|_2} - {\tau _i}{\text{ = }}{\left\| {\tilde z_j^i - \tilde x_j^i} \right\|_2} \le \left( {1/{\tau _i}} \right)\left\| {\tilde x_i^i - \tilde x_j^i} \right\|_2^2$ otherwise. We next bound the term $\left\| {\sum\nolimits_{j \in {\mathcal{B}_i}} {{w_{ij}}\left( {\tilde z_j^i - \tilde x_i^i} \right)} } \right\|_2^2$ in the following
\begin{equation}\label{AE7-1-3}
\left\| {\sum\limits_{j \in {\mathcal{B}_i}} {{w_{ij}}\left( {\tilde z_j^i - \tilde x_i^i} \right)} } \right\|_2^2 \le {\left( {\sum\limits_{j \in {\mathcal{B}_i}} {{w_{ij}}{\tau _i}} } \right)^2},
\end{equation}
where we use the fact that ${\left\| {\tilde z_j^i - \tilde x_i^i} \right\|_2} \le {\tau _i}$ if no clipping happens and ${\left\| {\tilde z_j^i - \tilde x_i^i} \right\|_2}{\text{ = }}{\tau _i}$ otherwise. To proceed, we fix the clipping parameter as ${\tau _i} = \sqrt {\left( {1/\sum\nolimits_{j \in {\mathcal{B}_i}} {{w_{ij}}} } \right)\sum\nolimits_{j \in {\mathcal{R}_i}} {{w_{ij}}\left\| {\tilde x_i^i - \tilde x_j^i} \right\|_2^2} } $ such that substituting (\ref{AE7-1-2}) and (\ref{AE7-1-3}) back into (\ref{AE7-1-1}) obtains
\begin{equation}\label{AE7-1-4}
\begin{aligned}
&\left\| {SC{C_i}\left( {\tilde x_i^i,{{\left\{ {\tilde x_j^i} \right\}}_{j \in {\mathcal{R}_i} \cup {\mathcal{B}_i}}}} \right) - {{\hat x}_i}} \right\|_2^2 \\
\le & 4\sum\limits_{b \in {\mathcal{B}_i}} {{w_{ib}}} \sum\limits_{r \in {\mathcal{R}_i}} {{w_{ir}}\left\| {\tilde x_i^i - \tilde x_r^i} \right\|_2^2}\\
\le & 16\sum\limits_{b \in {{\cal B}_i}} {{w_{ib}}\sum\limits_{r \in {{\cal R}_i}} {{w_{ir}}\mathop {\max }\limits_{j \in {{\cal R}_i} \cup \left\{ i \right\}} \left\| {\tilde x_j^i - {{\hat x}_i}} \right\|_2^2} } .
\end{aligned}
\end{equation}

The proof is completed via taking the square root on the both sides of (\ref{AE7-1-4}).

\subsection{Proof of Lemma \ref{L2}}\label{Section 7-2}
Define ${T_1}: = \nabla {f_i}\left( {{x_{i,k}},{\xi _{i,k}}} \right) - \left( {1/\left| \mathcal{R} \right|} \right)\sum\nolimits_{j \in \mathcal{R}} {\nabla {f_j}\left( {{x_{j,k}},{\xi _{j,k}}} \right)} $ and recall the definition of ${{\tilde D}_k}$ such that
\begin{equation*}
\begin{aligned}
\mathbb{E}{{\tilde D}_k}=& \sum\limits_{i \in \mathcal{R}} {\mathbb{E}\left\| {\tilde x_{i,k}^i - {{{\overset{\lower0.5em\hbox{$\smash{\scriptscriptstyle\frown}$}}{x} }_k}}} \right\|_2^2} \\
\le & \frac{1}{{1 - \eta }}\sum\limits_{i \in \mathcal{R}} {\mathbb{E}\left\| {{x_{i,k}} - {{\bar x}_k}} \right\|_2^2} + \frac{{2}}{\eta }\alpha _k^2\sum\limits_{i \in \mathcal{R}} {\mathbb{E}\left\| {{T_1}} \right\|_2^2}\\
& + \frac{{2}}{\eta }\alpha _k^2\sum\limits_{i \in \mathcal{R}} {\mathbb{E}\left\| {{{\tilde n}_{i,k}} - \sum\limits_{j \in \mathcal{R}} {\frac{{{{\tilde n}_{j,k}}}}{{\left| \mathcal{R} \right|}}} } \right\|_2^2}\\
\le & \frac{1}{{1 - \eta }}\sum\limits_{i \in \mathcal{R}} {\mathbb{E}\left\| {{x_{i,k}} - {{\bar x}_k}} \right\|_2^2}+ \frac{{2}}{\eta }\alpha _k^2\sum\limits_{i \in \mathcal{R}} {\mathbb{E}\left\| {{T_1}} \right\|_2^2}\\
&+ \frac{{2}}{\eta }\alpha _k^2\sum\limits_{i \in \mathcal{R}} {\mathbb{E}\left\langle {\left( {1 - \frac{1}{{\left| \mathcal{R} \right|}}} \right){{\tilde n}_{i,k}},\sum\limits_{j \in \mathcal{R}\backslash \left\{ i \right\}} {{{\tilde n}_{j,k}}} } \right\rangle } \\
& + \frac{2}{\eta }{\left( {1 - \frac{1}{{\left| \mathcal{R} \right|}}} \right)^2}\alpha _k^2\sum\limits_{i \in \mathcal{R}} {\mathbb{E}\left\| {{{\tilde n}_{i,k}}} \right\|_2^2} \\
\end{aligned}
\end{equation*}
\begin{equation}\label{AE7-2-1}
\begin{aligned}
& + \frac{{2}}{{\eta {{\left| \mathcal{R} \right|}^2}}}\alpha _k^2\sum\limits_{i \in \mathcal{R}} {\mathbb{E}\left\| {\sum\limits_{j \in \mathcal{R}\backslash \left\{ i \right\}} {{{\tilde n}_{j,k}}} } \right\|_2^2}\\
\le & \frac{1}{{1 - \eta }}\sum\limits_{i \in \mathcal{R}} \mathbb{E}D_k  + \frac{{2n\left| \mathcal{R} \right|}}{\eta }\varpi ^2 \alpha _k^2 + \frac{{2}}{\eta }\alpha _k^2\sum\limits_{i \in \mathcal{R}} {\mathbb{E}\left\| {{T_1}} \right\|_2^2},
\end{aligned}
\end{equation}
where the second inequality applies the update of Algorithm \ref{Algorithm 1} and the last inequality uses the fact that $\mathbb{E}{{\tilde n}_{i,k}} = 0$ and $\mathbb{E}\left\| {{{\tilde n}_{i,k}} - \mathbb{E}{{\tilde n}_{i,k}}} \right\|_2^2 = n\varpi ^2$. According to the standard variance decomposition,
\begin{equation}\label{AE7-2-2}
\mathbb{E}\left\| {{T_1}} \right\|_2^2 = \left\| {\mathbb{E}{T_1}} \right\|_2^2 + \mathbb{E}\left\| {{T_1} - \mathbb{E}{T_1}} \right\|_2^2.
\end{equation}

We next seek an upper bound on ${\left\| {\mathbb{E}{T_1}} \right\|_2^2}$ as follows:
\begin{equation}\label{AE7-2-3}
\begin{aligned}
&{\left\| {\mathbb{E}{T_1}} \right\|_2^2}\\
\le & 2\mathbb{E}\left\| {\nabla {f_i}\left( {{x_{i,k}}} \right) - \nabla {f_i}\left( {{{\bar x}_k}} \right)} \right\|_2^2 + 4\mathbb{E}\left\| {\nabla {f_i}\left( {{{\bar x}_k}} \right) - \nabla \bar F\left( {{{\bar x}_k}} \right)} \right\|_2^2 \\
&+ 4\mathbb{E}\left\| {\nabla \bar F\left( {{{\bar x}_k}} \right) - \frac{1}{{\left| \mathcal{R} \right|}}\sum\limits_{i \in \mathcal{R}} {\nabla {f_i}\left( {{x_{i,k}}} \right)} } \right\|_2^2\\
\le & 2{L^2}\mathbb{E}\left\| {{x_{i,k}} - {{\bar x}_k}} \right\|_2^2 + 4\mathbb{E}\left\| {\nabla {f_i}\left( {{{\bar x}_k}} \right) - \nabla \bar F\left( {{{\bar x}_k}} \right)} \right\|_2^2\\
& + \frac{{4{L^2}}}{{\left| \mathcal{R} \right|}}\mathbb{E}D_k\\
\le & \frac{{4{L^2}}}{{\left| \mathcal{R} \right|}}\mathbb{E}D_k + 2{L^2}\mathbb{E}\left\| {{x_{i,k}} - {{\bar x}_k}} \right\|_2^2 + 4{\zeta ^2},
\end{aligned}
\end{equation}
where the first inequality utilizes the basic inequality $\left\| {\tilde x + \tilde y} \right\|_2^2 \le 2\left\| {\tilde x} \right\|_2^2 + 2\left\| {\tilde y} \right\|_2^2$, $\forall \tilde x,\tilde y \in {\mathbb{R}^n}$ twice, the second inequality applies the $L$-smoothness (\ref{E2-1-2}) under Assumption \ref{A3}, and the last inequality is according to the bounded heterogeneity (\ref{E2-1-4}). We proceed to find an upper bound for $\mathbb{E}\left\| {{T_1} - \mathbb{E}{T_1}} \right\|_2^2$.
\begin{equation}\label{AE7-2-4}
\begin{aligned}
\mathbb{E}\left\| {{T_1} - \mathbb{E}{T_1}} \right\|_2^2 \le & 2\mathbb{E}\left\| {\sum\limits_{j \in \mathcal{R}} {\frac{{\nabla {f_j}\left( {{x_{j,k}},{\xi _{j,k}}} \right) - \nabla {f_j}\left( {{x_{j,k}}} \right)}}{{\left| \mathcal{R} \right|}}} } \right\|_2^2\\
&+ {\text{2}}\mathbb{E}\left\| {\nabla {f_i}\left( {{x_{i,k}},{\xi _{i,k}}} \right) - \nabla {f_i}\left( {{x_{i,k}}} \right)} \right\|_2^2\\
\le & {\text{4}}{\sigma ^2},
\end{aligned}
\end{equation}
where the first inequality utilizes the basic inequality and the last inequality is owing to the bounded variance (\ref{E2-1-3}). Combining (\ref{AE7-2-2}), (\ref{AE7-2-3}), and (\ref{AE7-2-4}) yields
\begin{equation}\label{AE7-2-5}
\sum\limits_{i \in \mathcal{R}} {\mathbb{E}\left\| {{T_1}} \right\|_2^2} =  6{L^2}\mathbb{E}D_k + 4\left| \mathcal{R} \right|\left( {{\sigma ^2} + {\zeta ^2}} \right).
\end{equation}

Plugging (\ref{AE7-2-5}) back into (\ref{AE7-2-1}) finishes the proof.

\subsection{Proof of Theorem \ref{T1}}\label{Section 7-3}
Recall the definition of ${D_{k + 1}} $ such that for any constant $\gamma  \in \left( {0,1} \right)$, we have
\begin{equation}\label{AE7-3-1}
\begin{aligned}
\mathbb{E}{D_{k + 1}} = & \mathbb{E}\left\| {\left( {{\mathbf{I}} - \frac{1}{{\left| \mathcal{R} \right|}}{\mathbf{1}}{{\mathbf{1}}^ \top }} \right)\left( {{X_{k + 1}} - \tilde W{{\tilde X}_k} + \tilde W{{\tilde X}_k}} \right)} \right\|_F^2 \\
\le &\frac{1}{{1 - \gamma }}\mathbb{E}\left\| {\tilde W{{\tilde X}_k} - \frac{1}{{\left| \mathcal{R} \right|}}{\mathbf{1}}{{\mathbf{1}}^ \top }\tilde W{{\tilde X}_k}} \right\|_F^2\\
& + \frac{1}{\gamma }\mathbb{E}\left\| {{X_{k + 1}} - {\tilde W}{{\tilde X}_k}} \right\|_F^2,
\end{aligned}
\end{equation}
where the inequality applies the following relations
\begin{equation}\label{AE7-3-2}
\left\| {{M_1} + {M_2}} \right\|_F^2 \le \frac{{\left\| {{M_1}} \right\|_F^2}}{{1 - \gamma }} + \frac{{\left\| {{M_2}} \right\|_F^2}}{\gamma },
\end{equation}
and the fact that $\left\| {{M_1}{M_2}} \right\|_F^2 \le \left\| {{M_1}} \right\|_2^2\left\| {{M_2}} \right\|_F^2$, for arbitrary matrices $M_1$ and $M_2$ with a same dimension, together with $\left\| {{\mathbf{I}} - \frac{1}{{\left| \mathcal{R} \right|}}{\mathbf{1}}{{\mathbf{1}}^ \top }} \right\|_2^2 = 1$. We proceed to bound $\mathbb{E}\left\| {{\tilde W}{{\tilde X}_k} - \frac{1}{{\left| \mathcal{R} \right|}}{\mathbf{1}}{{\mathbf{1}}^ \top }{\tilde W}{{\tilde X}_k}} \right\|_F^2$ in the sequel.
\begin{equation}\label{AE7-3-3}
\begin{aligned}
&\mathbb{E}\left\| {{\tilde W}{{\tilde X}_k} - \frac{1}{{\left| \mathcal{R} \right|}}{\mathbf{1}}{{\mathbf{1}}^ \top }{\tilde W}{{\tilde X}_k}} \right\|_F^2 \\
\le & \mathbb{E}\left\| {\left( {{\mathbf{I}} - \frac{1}{{\left| \mathcal{R} \right|}}{\mathbf{1}}{{\mathbf{1}}^ \top }} \right)\tilde W} \right\|_2^2\mathbb{E}\left\| {\left( {{\mathbf{I}} - \frac{1}{{\left| \mathcal{R} \right|}}{\mathbf{1}}{{\mathbf{1}}^ \top }} \right){{\tilde X}_k}} \right\|_F^2\\
= & \left( {1 - \lambda } \right){{\tilde D}_k},
\end{aligned}
\end{equation}
where the inequality applies the norm compatibility, i.e., for two arbitrary matrices $A \in {\mathbb{R}^{m \times n}}$ and $B \in {\mathbb{R}^{n \times d}}$, ${\left\| {AB} \right\|_F} \le {\left\| A \right\|_2}{\left\| B \right\|_F}$. According to \cite{Xin2022}, it can be verified that $0 < 1 - \lambda  \le 1$ under Assumption \ref{A1}. Considering the relation (\ref{E4-1-2}) in Definition \ref{D0}, we next seek an upper bound on $\mathbb{E}\left\| {{X_{k + 1}} - W{{\tilde X}_k}} \right\|_F^2$ as follows:
\begin{equation}\label{AE7-3-4}
\begin{aligned}
&\mathbb{E}\left\| {{X_{k + 1}} - W{{\tilde X}_k}} \right\|_F^2\\
= & \sum\limits_{i \in \mathcal{R}} {\mathbb{E}\left\| {{SCC_i}\left( {\tilde x_{i,k}^i,{{\left\{ {\tilde x_{j,k}^i} \right\}}_{j \in {\mathcal{R}_i} \cup {\mathcal{B}_i}}}} \right) - {{\hat x}_{i,k}}} \right\|_2^2} \\
\le & {\rho ^2}\sum\limits_{i \in \mathcal{R}} {\mathop {\max }\limits_{j \in {\mathcal{R}_i} \cup \left\{ i \right\}} } \mathbb{E}\left\| {\tilde x_{j,k}^i - {{\hat x}_{i,k}}} \right\|_2^2 \\
\le & 2{\rho ^2}\sum\limits_{i \in \mathcal{R}} {\mathop {\max }\limits_{j \in {\mathcal{R}_i} \cup \left\{ i \right\}} \mathbb{E}\left\| {\tilde x_{j,k}^j - {{{\overset{\lower0.5em\hbox{$\smash{\scriptscriptstyle\frown}$}}{x} }_k}}} \right\|_2^2}  + 2{\rho ^2}\sum\limits_{i \in \mathcal{R}} {\left\| {{{{\overset{\lower0.5em\hbox{$\smash{\scriptscriptstyle\frown}$}}{x} }_k}} - {{\hat x}_{i,k}}} \right\|_2^2} \\
\le & 2{\rho ^2}\sum\limits_{i \in \mathcal{R}} {\mathop {\max }\limits_{j \in \mathcal{R}} \mathbb{E}\left\| {\tilde x_{j,k}^j - {{{\overset{\lower0.5em\hbox{$\smash{\scriptscriptstyle\frown}$}}{x} }_k}}} \right\|_2^2}  \!+\! 2{\rho ^2}\sum\limits_{i \in \mathcal{R}} {\mathop {\max }\limits_{j \in \mathcal{R}} \mathbb{E}\left\| {\tilde x_{j,k}^j - {{{\overset{\lower0.5em\hbox{$\smash{\scriptscriptstyle\frown}$}}{x} }_k}}} \right\|_2^2} \\
\le & 4\left| \mathcal{R} \right|{\rho ^2}\sum\limits_{i \in \mathcal{R}} {\mathbb{E}\left\| {\tilde x_{i,k}^i - {{{\overset{\lower0.5em\hbox{$\smash{\scriptscriptstyle\frown}$}}{x} }_k}}} \right\|_2^2} \\
= & 4\left| \mathcal{R} \right|{\rho ^2}\mathbb{E}{{\tilde D}_k}.
\end{aligned}
\end{equation}

Substituting (\ref{AE7-3-3}) and (\ref{AE7-3-4}) back into (\ref{AE7-3-1}) yields
\begin{equation}\label{AE7-3-5}
\mathbb{E}{D_{k + 1}}\le \left( {\frac{{1 - \lambda }}{{1 - \gamma }} + \frac{{4\left| \mathcal{R} \right|}}{\gamma }{\rho ^2}} \right)\mathbb{E}{{\tilde D}_k}.
\end{equation}

We choose $0 \le \rho  < \lambda /\left( {4\sqrt {\left| \mathcal{R} \right|} } \right)$ and let $\gamma = 2\rho \sqrt {\left| \mathcal{R} \right|} $ such that combining (\ref{E4-1-3}) and (\ref{AE7-3-5}) yields
\begin{equation}\label{AE7-3-6}
\begin{aligned}
\mathbb{E}{D_{k + 1}} \le & \left( {1 + 2\gamma  - \lambda } \right)\left( {\frac{1}{{1 - \eta }} + \frac{{12{L^2}}}{\eta }\alpha _k^2} \right)\mathbb{E}{D_k}\\
&+ \left( {1 + 2\gamma  - \lambda } \right)\frac{{2\left| \mathcal{R} \right|}}{\eta }\left( {n\varpi ^2 + 4\left( {{\sigma ^2} + {\zeta ^2}} \right)} \right)\alpha _k^2\\
\le &\left( {1 - \varphi } \right)\frac{{2\left| \mathcal{R} \right|}}{\eta }\left( {n\varpi ^2 + 4\left( {{\sigma ^2} + {\zeta ^2}} \right)} \right)\alpha _k^2\\
& + \left( {1 - \varphi } \right)\left( {\frac{1}{{1 - \eta }} { + \frac{{12{L^2}}}{\eta }\alpha _k^2} } \right)\mathbb{E}{D_k}.
\end{aligned}
\end{equation}

If we further fix $\eta  = \varphi /2$ and choose the step-size $0 < {\alpha _k} \le \varphi /\left( {4L\sqrt 3 \left( {4 - \varphi } \right)} \right)$, then (\ref{AE7-3-6}) becomes
\begin{equation}\label{AE7-3-7}
\mathbb{E}{D_{k + 1}} \le  {\frac{\varphi \mathbb{E}{D_k} }{{4 - \varphi }}} + \frac{{4\left( {1 - \varphi } \right)\left| \mathcal{R} \right|}}{\varphi }\left( {n\varpi ^2 + 4\left( {{\sigma ^2} + {\zeta ^2}} \right)} \right)\alpha _k^2.
\end{equation}

Via defining $\phi := \varphi /\left( {4 - \varphi } \right)$ and $\vartheta := 4\left| \mathcal{R} \right|\left( {1 - \varphi } \right)\left( {n{\varpi ^2} + 4\left( {{\sigma ^2} + {\zeta ^2}} \right)} \right)/\varphi $, (\ref{AE7-3-7}) reduces to
\begin{equation}\label{AE7-3-8}
\mathbb{E}{D_{k + 1}} \le \left( {1 - \phi} \right)\mathbb{E}{D_k} + \vartheta \alpha _k^2.
\end{equation}

If we choose a decaying step-size ${\alpha _k} = \theta /\left( {k + {k_0}} \right)$, then applying telescopic cancellation on (\ref{AE7-3-8}) obtains
\begin{equation}\label{AE7-3-9}
\begin{aligned}
\mathbb{E}{D_k} \le &  {\left( {1 - \phi } \right)^k}\mathbb{E}{D_0} + \frac{{\vartheta {\theta ^2}}}{{{{\left( {k + {k_0} - 1} \right)}^2}}} + \frac{{\left( {1 - \phi } \right)\vartheta {\theta ^2}}}{{{{\left( {k + {k_0} - 2} \right)}^2}}}\\
& +  \cdots  + \frac{{\vartheta {\theta ^2}{{\left( {1 - \phi } \right)}^{k - 1}}}}{{k_0^2}}.
\end{aligned}
\end{equation}

According to \cite[Lemma 5]{Wu2023}, there exists a constant $\iota $ satisfying $\iota  \ge{\left( {{k_0} + 1} \right)^2}/k_0^2$ such that
\begin{equation}\label{AE7-3-10}
\mathbb{E}{D_k} \le {\left( {1 - \phi } \right)^k}{D_0} + \frac{{2\iota \vartheta {\theta ^2}}}{\phi }\frac{1}{{{{\left( {k + {k_0}} \right)}^2}}},
\end{equation}
which is exactly the first result (\ref{E4-1-4}). We then fix the step-size ${\alpha _k} \equiv \alpha $ and update (\ref{AE7-3-8}) recursively to get
\begin{equation}\label{AE7-3-11}
\begin{aligned}
\mathbb{E}{D_{k + 1}} \le & \left( {1 - \phi } \right)\mathbb{E}{D_k} + \vartheta {\alpha ^2} \\
\le & {\left( {1 - \phi } \right)^{k + 1}}{D_0} + \vartheta {\alpha ^2}\sum\limits_{t = 0}^k {{{\left( {1 - \phi } \right)}^k}} \\
\le & {\left( {1 - \phi } \right)^{k + 1}}{D_0} + \frac{\vartheta }{\phi }{\alpha ^2},
\end{aligned}
\end{equation}
which verifies the second result (\ref{E4-1-5}). This completes the proof.

\subsection{Proof of Theorem \ref{T2}}\label{Section 7-4}
Under Assumption \ref{A3}, we know that the global objective is $L$-smooth such that
\begin{equation}\label{AE7-4-1}
\begin{aligned}
\mathbb{E}f\left( {{{\bar x}_{k + 1}}} \right) \le & \mathbb{E}f\left( {{{\bar x}_k}} \right) + \mathbb{E}\left\langle {\nabla f\left( {{{\bar x}_k}} \right),{{\bar x}_{k + 1}} - {{\bar x}_k}} \right\rangle \\
&+ \frac{L}{2}\mathbb{E}\left\| {{{\bar x}_{k + 1}} - {{\bar x}_k}} \right\|_2^2,
\end{aligned}
\end{equation}
which is also known as the standard descent lemma. We next seek an upper bound for $\mathbb{E}\left\| {{{\bar x}_{k + 1}} - {{\bar x}_k}} \right\|_2^2$ in the right-hand-side (RHS) of (\ref{AE7-4-1}) as follows:
\begin{equation}\label{AE7-4-2}
\begin{aligned}
& \mathbb{E}\left\| {{{\bar x}_{k + 1}} - {{\bar x}_k}} \right\|_2^2\\
= & \left\| {{{\bar x}_{k + 1}} - {{\bar x}_k} + {\alpha _k}} \right.\left( {\nabla f\left( {{{\bar x}_k};{\xi _k}} \right) - \nabla f\left( {{{\bar x}_k};{\xi _k}} \right)} \right)\\
&\left. { + {\alpha _k}\left( {\nabla f\left( {{{\bar x}_k}} \right) - \nabla f\left( {{{\bar x}_k}} \right)} \right)} \right\|_2^2\\
\le & 2\alpha _k^2\mathbb{E}\left\| {\frac{1}{{{\alpha _k}}}\left( {{{\bar x}_{k + 1}} - {{\bar x}_k}} \right) + \nabla f\left( {{{\bar x}_k};{\xi _k}} \right) - \nabla f\left( {{{\bar x}_k}} \right)} \right\|_2^2\\
&+ 2\alpha _k^2\left\| {\nabla f\left( {{{\bar x}_k};{\xi _k}} \right) - \nabla f\left( {{{\bar x}_k}} \right)} \right\|_2^2\\
\le & 2\alpha _k^2\left\| {\frac{1}{{{\alpha _k}}}\left( {{{\bar x}_{k + 1}} - {{\bar x}_k}} \right) + \nabla f\left( {{{\bar x}_k};{\xi _k}} \right) - \nabla f\left( {{{\bar x}_k}} \right)} \right\|_2^2\\
& + \frac{2}{{\left| \mathcal{R} \right|}}\alpha _k^2\sum\limits_{i \in \mathcal{R}} {\mathbb{E}\left\| {\nabla {f_i}\left( {{{\bar x}_k};{\xi _k}} \right) - \nabla {f_i}\left( {{{\bar x}_k}} \right)} \right\|_2^2} \\
\le & 2\alpha _k^2\mathbb{E}\left\| {\frac{1}{{{\alpha _k}}}\left( {{{\bar x}_{k + 1}} - {{\bar x}_k}} \right) + \nabla f\left( {{{\bar x}_k};{\xi _k}} \right) - \nabla f\left( {{{\bar x}_k}} \right)} \right\|_2^2\\
& + 2{\sigma ^2}\alpha _k^2,
\end{aligned}
\end{equation}
where the second inequality applies the basic inequality, the third inequality leverages the Jensen's inequality, and the last inequality employs the bounded variance (\ref{E2-1-3}) under Assumption \ref{A5}-i). We proceed to bound $\mathbb{E}\left\langle {\nabla f\left( {{{\bar x}_k}} \right),{{\bar x}_{k + 1}} - {{\bar x}_k}} \right\rangle $ in the RHS of (\ref{AE7-4-1}).
\begin{equation}\label{AE7-4-3}
\begin{aligned}
&\mathbb{E}\left\langle {\nabla f\left( {{{\bar x}_k}} \right),{{\bar x}_{k + 1}} - {{\bar x}_k}} \right\rangle \\
= & {\alpha _k}\mathbb{E}\left\langle {\nabla f\left( {{{\bar x}_k}} \right),\nabla f\left( {{{\bar x}_k};{\xi _k}} \right) \!-\! \nabla f\left( {{{\bar x}_k}} \right) \!+\! \frac{1}{{{\alpha _k}}}\left( {{{\bar x}_{k + 1}} - {{\bar x}_k}} \right)} \right\rangle \\
= & \frac{{{\alpha _k}}}{2}\mathbb{E}\left\| {\nabla f\left( {{{\bar x}_k};{\xi _k}} \right) + \frac{1}{{{\alpha _k}}}\left( {{{\bar x}_{k + 1}} \!-\! {{\bar x}_k}} \right)} \right\|_2^2 \!-\! \frac{{{\alpha _k}}}{2}\mathbb{E}\left\| {\nabla f\left( {{{\bar x}_k}} \right)} \right\|_2^2\\
&- \frac{{{\alpha _k}}}{2}\mathbb{E}\left\| {\nabla f\left( {{{\bar x}_k};{\xi _k}} \right) - \nabla f\left( {{{\bar x}_k}} \right) + \frac{1}{{{\alpha _k}}}\left( {{{\bar x}_{k + 1}} - {{\bar x}_k}} \right)} \right\|_2^2,
\end{aligned}
\end{equation}
where the first equality applies the fact that $\mathbb{E}\left\langle {\nabla f\left( {{{\bar x}_k}} \right),\nabla f\left( {{{\bar x}_k};{\xi _k}} \right) - \nabla f\left( {{{\bar x}_k}} \right)} \right\rangle  = 0$ and the second equality follows the Cosine theorem, i.e.,
$\left\langle {\tilde x, \tilde y} \right\rangle  = \frac{1}{2}\left\| {\tilde x + \tilde y} \right\|_2^2 - \frac{1}{2}\left\| \tilde x \right\|_2^2 - \frac{1}{2}\left\| \tilde y \right\|_2^2$, $\forall \tilde x,\tilde y \in {\mathbb{R}^n}$. We next substitute (\ref{AE7-4-2}) and (\ref{AE7-4-3}) back into (\ref{AE7-4-1}) to obtain
\begin{equation}\label{AE7-4-4}
\begin{aligned}
\mathbb{E}f\left( {{{\bar x}_{k + 1}}} \right)\le & \mathbb{E}f\left( {{{\bar x}_k}} \right) + \frac{{{\alpha _k}}}{2}\mathbb{E}\left\| {\nabla f\left( {{{\bar x}_k};{\xi _k}} \right) + \frac{{{{\bar x}_{k + 1}} - {{\bar x}_k}}}{{{\alpha _k}}}} \right\|_2^2\\
& - \frac{{{\alpha _k}}}{2}\mathbb{E}\left\| {\nabla f\left( {{{\bar x}_k}} \right)} \right\|_2^2 + L{\sigma ^2}\alpha _k^2.
\end{aligned}
\end{equation}

We continue to define ${V_1}: = \nabla f\left( {{{\bar x}_k};{\xi _k}} \right) - \left( {1/\left| \mathcal{R} \right|} \right) \sum\nolimits_{j \in \mathcal{R}} {\nabla {f_j}\left( {{x_{j,k}};{\xi _{j,k}}} \right)} $, ${V_2}: = \left( {1/\left( {\left| \mathcal{R} \right|{\alpha _k}} \right)} \right)$ $\sum\nolimits_{i \in \mathcal{R}} {\left( {{{\hat x}_{i,k}} - {{\bar x}_k} + \left( {{\alpha _k}/\left| \mathcal{R} \right|} \right)\sum\nolimits_{j \in \mathcal{R}} {\nabla {f_j}\left( {{x_{j,k}};{\xi _{j,k}}} \right)} } \right)} $, and ${V_3}: = \left( {1/\left( {\left| \mathcal{R} \right|{\alpha _k}} \right)} \right)\sum\nolimits_{i \in \mathcal{R}} {\left( {SC{C_i}\left\{ {\tilde x_{i,k}^i,{{\left\{ {\tilde x_{j,k}^i} \right\}}_{j \in {\mathcal{R}_i} \cup \left\{ i \right\}}}} \right\}} \right.}  $ $\left. { - {{\hat x}_{i,k}}} \right)$. According to the update rule of Algorithm \ref{Algorithm 1}, we expand $\nabla f\left( {{{\bar x}_k};{\xi _k}} \right) + \left( {{{\bar x}_{k + 1}} - {{\bar x}_k}} \right)/{\alpha _k}$ in the RHS of (\ref{AE7-4-3}) as follows:
\begin{equation}\label{AE7-4-5}
\begin{aligned}
&\nabla f\left( {{{\bar x}_k};{\xi _k}} \right) + \frac{1}{{{\alpha _k}}}\left( {{{\bar x}_{k + 1}} - {{\bar x}_k}} \right) \\
= & \frac{1}{{\left| \mathcal{R} \right|{\alpha _k}}}\sum\limits_{i \in \mathcal{R}} {\left( {SC{C_i}\left( {\tilde x_{i,k}^i,{{\left\{ {\tilde x_{j,k}^i} \right\}}_{j \in {\mathcal{R}_i} \cup \left\{ i \right\}}}} \right) - {{\bar x}_k}} \right)}\\
&+\nabla f\left( {{{\bar x}_k};{\xi _k}} \right)\\
= & V_1 + V_2 + V_3.
\end{aligned}
\end{equation}

We next seek an upper bound for $\mathbb{E}\left\| {{V_1}} \right\|_2^2$ as follows:
\begin{equation}\label{AE7-4-6}
\begin{aligned}
\mathbb{E}\left\| {{V_1}} \right\|_2^2 = & \mathbb{E}\left\| {\frac{1}{{\left| \mathcal{R} \right|}}\sum\limits_{i \in \mathcal{R}} {\left( {\nabla {f_i}\left( {{{\bar x}_k};{\xi _k}} \right) - \nabla {f_i}\left( {{x_{i,k}};{\xi _{i,k}}} \right)} \right)} } \right\|_2^2 \\
\le & \frac{1}{{\left| \mathcal{R} \right|}}\sum\limits_{i \in \mathcal{R}} {\mathbb{E}\left\| {\nabla {f_i}\left( {{{\bar x}_k};{\xi _k}} \right) - \nabla {f_i}\left( {{x_{i,k}};{\xi _{i,k}}} \right)} \right\|_2^2}\\
\le & \frac{{{L^2}}}{{\left| \mathcal{R} \right|}}\sum\limits_{i \in \mathcal{R}} {\mathbb{E}\left\| {{x_{i,k}} - {{\bar x}_k}} \right\|_2^2} \\
= & \frac{{{L^2}}}{{\left| {\cal R} \right|}}\left\| {{X_k} - \frac{1}{{\left| {\cal R} \right|}}{\bf{1}}{{\bf{1}}^ \top }{X_k}} \right\|_2^2\\
= & \frac{{{L^2}}}{{\left| \mathcal{R} \right|}}\mathbb{E}{D_k},
\end{aligned}
\end{equation}
where the second and third inequalities apply the Jensen's inequality and the $L$-smoothness (\ref{E2-1-2}) under Assumption \ref{A3}, respectively. According to the algorithm update (\ref{E3-3}), we next bound $\mathbb{E}\left\| {{V_2}} \right\|_2^2$ as follows:
\begin{equation}\label{AE7-4-7}
\begin{aligned}
\mathbb{E}\left\| {{V_2}} \right\|_2^2 = & \frac{1}{{{{\left| \mathcal{R} \right|}^2}\alpha _k^2}}\mathbb{E}\left\| {\sum\limits_{i \in \mathcal{R}} {\left( {{{\hat x}_{i,k}} - {{{\overset{\lower0.5em\hbox{$\smash{\scriptscriptstyle\frown}$}}{x} }_k}}} \right)}  - {\alpha _k}\sum\limits_{j \in \mathcal{R}} {{{\tilde n}_{j,k}}} } \right\|_2^2\\
\le &\frac{2}{{{{\left| \mathcal{R} \right|}^2}\alpha _k^2}}\mathbb{E}\left\| {\sum\limits_{i \in \mathcal{R}} {\left( {{{\hat x}_{i,k}} - {{{\overset{\lower0.5em\hbox{$\smash{\scriptscriptstyle\frown}$}}{x} }_k}}} \right)} } \right\|_2^2 \!+\! \frac{2}{{{{\left| \mathcal{R} \right|}^2}}}\mathbb{E}\left\| {\sum\limits_{j \in \mathcal{R}} {{{\tilde n}_{j,k}}} } \right\|_2^2\\
\le & \frac{2}{{{{\left| \mathcal{R} \right|}^2}\alpha _k^2}}\mathbb{E}\left\| {\sum\limits_{i \in \mathcal{R}} {\left( {{{\hat x}_{i,k}} - {{{\overset{\lower0.5em\hbox{$\smash{\scriptscriptstyle\frown}$}}{x} }_k}}} \right)} } \right\|_2^2 + \frac{2}{{\left| \mathcal{R} \right|}}\sum\limits_{i \in \mathcal{R}} {\mathbb{E}\left\| {{{\tilde n}_{i,k}}} \right\|_2^2} \\
= & \frac{2}{{{{\left| \mathcal{R} \right|}^2}\alpha _k^2}}\mathbb{E}\left\| {\left( {{{\mathbf{1}}^ \top }\tilde W - {{\mathbf{1}}^ \top }} \right)\left( {{{\tilde X}_k} - \frac{1}{{\left| \mathcal{R} \right|}}{\mathbf{1}}{{\mathbf{1}}^ \top }{{\tilde X}_k}} \right)} \right\|_F^2\\
& + \frac{2}{{\left| \mathcal{R} \right|}}\sum\limits_{i \in \mathcal{R}} {\mathbb{E}\left\| {{{\tilde n}_{i,k}}} \right\|_2^2}\\
\leq & \frac{2}{{{{\left| \mathcal{R} \right|}^2}\alpha _k^2}}\mathbb{E}\left\| {{{\mathbf{1}}^ \top }\tilde W - {{\mathbf{1}}^ \top }} \right\|_2^2\mathbb{E}{{\tilde D}_k} \!+\! \frac{2}{{\left| \mathcal{R} \right|}}\sum\limits_{i \in \mathcal{R}} {\mathbb{E}\left\| {{{\tilde n}_{i,k}}} \right\|_2^2} \\
\le & 2n\varpi ^2,
\end{aligned}
\end{equation}
where the second inequality applies the basic inequality, the second inequality is owing to the Jensen's inequality, the third inequality follows the norm compatibility again, and the last inequality uses the fact that ${\left\| {{{\mathbf{1}}^ \top }\tilde W - {{\mathbf{1}}^ \top }} \right\|_2} = 0$ since ${\tilde W}$ is doubly stochastic according to (\ref{E4-1-1}). To proceed, an upper bound on the term $\mathbb{E}\left\| {{V_3}} \right\|_2^2 $ is sought in the following
\begin{equation}\label{AE7-4-8}
\begin{aligned}
\mathbb{E}\left\| {{V_3}} \right\|_2^2 \le & \frac{{{\rho ^2}}}{{\left| \mathcal{R} \right|\alpha _k^2}}\sum\limits_{i \in \mathcal{R}} {\mathop {\max }\limits_{j \in {\mathcal{R}_i} \cup \left\{ i \right\}} } \mathbb{E}\left\| {\tilde x_{j,k}^i - {{\hat x}_{i,k}}} \right\|_2^2 \\
\le & \frac{{2{\rho ^2}}}{{\left| \mathcal{R} \right|\alpha _k^2}}\sum\limits_{i \in \mathcal{R}} {\mathop {\max }\limits_{j \in {\mathcal{R}_i} \cup \left\{ i \right\}} \mathbb{E}\left\| {\tilde x_{j,k}^j - {{{\overset{\lower0.5em\hbox{$\smash{\scriptscriptstyle\frown}$}}{x} }_k}}} \right\|_2^2} \\
& + \frac{{2{\rho ^2}}}{{\alpha _k^2}}\mathop {\max }\limits_{j \in \mathcal{R}} \mathbb{E}\left\| {\tilde x_{j,k}^j - {{{\overset{\lower0.5em\hbox{$\smash{\scriptscriptstyle\frown}$}}{x} }_k}}} \right\|_2^2 \\
\le & \frac{{4{\rho ^2}}}{{\alpha _k^2}}\mathop {\max }\limits_{j \in \mathcal{R}} \mathbb{E}\left\| {\tilde x_{j,k}^j - {{{\overset{\lower0.5em\hbox{$\smash{\scriptscriptstyle\frown}$}}{x} }_k}}} \right\|_2^2 \\
\le & \frac{{4{\rho ^2}}}{{\alpha _k^2}}\sum\limits_{i \in \mathcal{R}} {\mathbb{E}\left\| {\tilde x_{i,k}^i - {{{\overset{\lower0.5em\hbox{$\smash{\scriptscriptstyle\frown}$}}{x} }_k}}} \right\|_2^2} \\
= & \frac{{4{\rho ^2}}}{{\alpha _k^2}}\mathbb{E}{{\tilde D}_k},
\end{aligned}
\end{equation}
where the first inequality utilizes the relation (\ref{E4-1-2}) in Definition \ref{D0} and the second inequality uses the basic inequality. To recap, plugging the relations (\ref{AE7-4-6})-(\ref{AE7-4-8}) back into (\ref{AE7-4-5}) yields
\begin{equation}\label{AE7-4-9}
\begin{aligned}
& \mathbb{E}\left\| {\nabla f\left( {{{\bar x}_k};{\xi _k}} \right) + \frac{1}{{{\alpha _k}}}\left( {{{\bar x}_{k + 1}} - {{\bar x}_k}} \right)} \right\|_2^2 \\
\le & 2\mathbb{E}\left\| {{V_1}} \right\|_2^2 + 4\mathbb{E}\left\| {{V_2}} \right\|_2^2 + 4\mathbb{E}\left\| {{V_3}} \right\|_2^2 \\
\le & \frac{{2{L^2}}}{{\left| \mathcal{R} \right|}}\mathbb{E}{D_k} + 8n\varpi ^2 + \frac{{16{\rho ^2}}}{{\alpha _k^2}}\mathbb{E}{{\tilde D}_k},
\end{aligned}
\end{equation}
where the first inequality applies the basic inequality twice. Plugging (\ref{E4-1-3}) into (\ref{AE7-4-9}) obtains
\begin{equation*}
\begin{aligned}
& \mathbb{E}\left\| {\nabla f\left( {{{\bar x}_k};{\xi _k}} \right) + \frac{1}{{{\alpha _k}}}\left( {{{\bar x}_{k + 1}} - {{\bar x}_k}} \right)} \right\|_2^2 \\
\end{aligned}
\end{equation*}
\begin{equation}\label{AE7-4-10}
\begin{aligned}
\le & 2\left( {\frac{{96\left| \mathcal{R} \right|{L^2}{\rho ^2}}}{\eta } + \frac{{{L^2}}}{{\left| \mathcal{R} \right|}} + \frac{{8{\rho ^2}}}{{1 - \eta }}\frac{1}{{\alpha _k^2}}} \right)\mathbb{E}{D_k} \\
& + \frac{{128\left| \mathcal{R} \right|}}{\eta }{\rho ^2}\left( {{\sigma ^2} + {\zeta ^2}} \right) + 8n\left( {1 + \frac{{8\left| \mathcal{R} \right|{\rho ^2}}}{\eta }} \right)\varpi ^2.
\end{aligned}
\end{equation}

We then substitute (\ref{AE7-4-10}) into (\ref{AE7-4-4}) to get
\begin{equation}\label{AE7-4-11}
\begin{aligned}
& \mathbb{E}f\left( {{{\bar x}_{k + 1}}} \right)\\
\le &  \mathbb{E}f\left( {{{\bar x}_k}} \right) + \left( {\frac{{96\left| \mathcal{R} \right|{L^2}{\rho ^2}}}{\eta } + \frac{{{L^2}}}{{\left| \mathcal{R} \right|}} + \frac{{8{\rho ^2}}}{{1 - \eta }}\frac{1}{{\alpha _k^2}}} \right){\alpha _k}\mathbb{E}{D_k}\\
&+4\left( {\frac{{16\left| \mathcal{R} \right| {\rho ^2}}}{\eta }\left( {{\sigma ^2} + {\zeta ^2}} \right) + n\left( {1 + \frac{{8\left| \mathcal{R} \right|{\rho ^2}}}{\eta }} \right){\varpi ^2}} \right){\alpha _k}\\
&- \frac{{{\alpha _k}}}{2}\mathbb{E}\left\| {\nabla f\left( {{{\bar x}_k}} \right)} \right\|_2^2 + L{\sigma ^2}\alpha _k^2.
\end{aligned}
\end{equation}

Applying the P-Ł condition (\ref{E2-1-5}) to the inequality (\ref{AE7-4-11}) becomes
\begin{equation}\label{AE7-4-12}
\begin{aligned}
& \mathbb{E}f\left( {{{\bar x}_{k + 1}}} \right) - {f^*}\\
\le & 4\left( {\frac{{16\left| \mathcal{R} \right| {\rho ^2}}}{\eta }\left( {{\sigma ^2} + {\zeta ^2}} \right)+n\left( {1 + \frac{{8\left| \mathcal{R} \right|{\rho ^2}}}{\eta }} \right){\varpi ^2}} \right){\alpha _k}\\
& + \left( {\frac{{96\left| \mathcal{R} \right|{L^2}{\rho ^2}}}{\eta } + \frac{{{L^2}}}{{\left| \mathcal{R} \right|}} + \frac{{8{\rho ^2}}}{{1 - \eta }}\frac{1}{{\alpha _k^2}}} \right){\alpha _k}\mathbb{E}{D_k}\\
&+\left( {1 - \nu \alpha_k } \right)\left( {\mathbb{E}f\left( {{{\bar x}_k}} \right) - {f^*}} \right) + L{\sigma ^2}\alpha _k^2.
\end{aligned}
\end{equation}

If we further choose the decaying step-size ${\alpha _k} = \underline{\theta} /\left( {k + {k_0}} \right)$ with $\underline{\theta} = \min \left\{ {1/\nu ,\phi /\left( {4\sqrt 3 L} \right)} \right\}$, summing (\ref{AE7-4-12}) over $k$ from 0 to $K$, $\forall K \ge 1$, yields
\begin{equation}\label{AE7-4-13}
\begin{aligned}
&\nu \sum\limits_{k = 0}^K {{\alpha _k}\left( {\mathbb{E}f\left( {{{\bar x}_k}} \right) - {f^*}} \right)} \\
\le & \left( {\frac{{96\left| \mathcal{R} \right|{\rho ^2}}}{\eta } + \frac{1}{{\left| \mathcal{R} \right|}}} \right){L^2}\sum\limits_{k = 0}^K {{\alpha _k}\mathbb{E}{D_k}} + \frac{{8{\rho ^2}}}{{1 - \eta }}\sum\limits_{k = 0}^K {\frac{1}{{{\alpha _k}}}\mathbb{E}{D_k}} \\
& + L{\sigma ^2}\sum\limits_{k = 0}^K {\alpha _k^2} + \mathbb{E}f\left( {{{\bar x}_0}} \right) - {f^*} - \left( {\mathbb{E}f\left( {{{\bar x}_{K + 1}}} \right) - {f^*}} \right) \\
& + 4\left( {\frac{{16\left| \mathcal{R} \right|{\rho ^2}}}{\eta }\left( {{\sigma ^2} + {\zeta ^2}} \right) \!+\! n\left( {1 \! + \! \frac{{8\left| \mathcal{R} \right|{\rho ^2}}}{\eta }} \right){\varpi ^2}} \right)\sum\limits_{k = 0}^K {{\alpha _k}}.
\end{aligned}
\end{equation}

Since $0 < \nu {\alpha _k} < 1$, we let $\mathbb{E}f_{K + 1}^{{\text{best}}} = {\min _{t \in \left\{ {1,2, \ldots ,K + 1} \right\}}}f\left( {{{\bar x}_t}} \right)$ such that $\mathbb{E}f_{K + 1}^{{\text{best}}} - {f^*} \ge 0$. We rearrange (\ref{AE7-4-13}) to generate
\begin{equation}\label{AE7-4-14}
\begin{aligned}
& {\mathbb{E}f_{K + 1}^{{\text{best}}} - {f^*}}\\
\le & \frac{{\mathbb{E}f\left( {{{\bar x}_0}} \right) - {f^*}}}{{\underline{\theta } \nu \left( {\ln\left( {K + {k_0}} \right) - \ln \left( {{k_0}} \right)} \right)}} + \frac{{\underline{\theta } L{\sigma ^2}\sum\limits_{k = 0}^K {\frac{1}{{{{\left( {k + {k_0}} \right)}^2}}}} }}{{\nu \left( {\ln \left( {K + {k_0}} \right) - \ln \left( {{k_0}} \right)} \right)}} \\
& + \frac{{{L^2}}}{\nu }\left( {\frac{{96\left| \mathcal{R} \right|{\rho ^2}}}{\eta } + \frac{1}{{\left| \mathcal{R} \right|}}} \right) \frac{{\sum\limits_{k = 0}^K {\frac{1}{{k + {k_0}}}\mathbb{E}{D_k}} }}{{\ln \left( {K + {k_0}} \right) - \ln \left( {{k_0}} \right)}}\\
& + \frac{{8{\rho ^2}}}{{\nu \left( {1 - \eta } \right){{\underline{\theta } }^2}}}\frac{{\sum\limits_{k = 0}^K {\left( {k + {k_0}} \right)\mathbb{E}{D_k}} }}{{\ln \left( {K + {k_0}} \right) - \ln \left( {{k_0}} \right)}} \!+\! \frac{{64\left| \mathcal{R} \right|{\rho ^2}}}{{\nu \eta }}\left( {{\sigma ^2} + {\zeta ^2}} \right)\\
& + \frac{{4n}}{\nu }\left( {1 + \frac{{8\left| \mathcal{R} \right|{\rho ^2}}}{\eta }} \right){\varpi ^2}.
\end{aligned}
\end{equation}
If $K$ approaches to infinity, then it follows from the relation (\ref{E4-1-4}) that (\ref{AE7-4-14}) gives rise to an asymptotic convergence error of $\mathcal{O}\left( {{\rho ^2}\left( {{\varpi ^2} + {\sigma ^2} + {\zeta ^2} } \right)} + {\varpi ^2} \right)$,
%\begin{equation}\label{AE7-4-15}
%\mathop {{\text{lim}}}\limits_{K \to \infty } \mathbb{E}f_{K + 1}^{{\text{best}}} - {f^*} \le \mathcal{O}\left( {{\rho ^2}\left( {{\sigma ^2} + {\zeta ^2} + {\varpi ^2}} \right)} \right),
%\end{equation}
which completes the proof.

\subsection{Proof of Corollary \ref{C4-3-1}}\label{Section 7-5}
Since $\varpi = \rho = 0$, it follows from (\ref{AE7-4-14}) that
\begin{equation}\label{AE7-5-1}
\begin{aligned}
& {\mathbb{E}f_{K + 1}^{{\text{best}}} - {f^*}}\\
\le &  \frac{1}{{\left| \mathcal{R} \right|}}\frac{{\sum\limits_{k = 0}^K {\frac{1}{{k + {k_0}}}\mathbb{E}{D_k}} }}{{\ln \left( {K + {k_0}} \right) - \ln \left( {{k_0}} \right)}}{\text{ }} + \frac{{\underline{\theta } L{\sigma ^2}\sum\limits_{k = 0}^K {\frac{1}{{{{\left( {k + {k_0}} \right)}^2}}}} }}{{\ln \left( {K + {k_0}} \right) - \ln \left( {{k_0}} \right)}}\\
&+\frac{{\mathbb{E}f\left( {{{\bar x}_0}} \right) - {f^*}}}{{\underline{\theta } \left( {\ln \left( {K + {k_0}} \right) - \ln \left( {{k_0}} \right)} \right)}}.
\end{aligned}
\end{equation}

In view of the relation (\ref{E4-1-4}) in Theorem \ref{T1}, the proof is completed via taking $K$ to infinity.

\subsection{Proof of Theorem \ref{T3}}\label{Section 7-6}
Following the same technical line as (\ref{AE7-4-1})-(\ref{AE7-4-11}), we set $\alpha_k \equiv \alpha$ such that (\ref{AE7-4-12}) becomes
\begin{equation}\label{AE7-6-17}
\begin{aligned}
& \mathbb{E}f\left( {{{\bar x}_{k + 1}}} \right) - {f^*}\\
\le & \left( {\frac{{96\left| \mathcal{R} \right|{L^2}{\rho ^2}}}{\eta } + \frac{{{L^2}}}{{\left| \mathcal{R} \right|}} + \frac{{8{\rho ^2}}}{{1 - \eta }}\frac{1}{{\alpha^2}}} \right){\alpha}\mathbb{E}{D_k} + L{\sigma ^2}\alpha^2\\
& + 4\left( {\frac{{16\left| \mathcal{R} \right| {\rho ^2}}}{\eta }\left( {{\sigma ^2} + {\zeta ^2}} \right)+n\left( {1 + \frac{{8\left| \mathcal{R} \right|{\rho ^2}}}{\eta }} \right){\varpi ^2}} \right){\alpha}\\
&+\left( {1 - \nu \alpha } \right)\left( {\mathbb{E}f\left( {{{\bar x}_k}} \right) - {f^*}} \right).
\end{aligned}
\end{equation}

We then rearrange (\ref{AE7-6-17}) to obtain
\begin{equation}\label{AE7-6-18}
\begin{aligned}
& \mathbb{E}f\left( {{{\bar x}_k}} \right) - {f^*}\\
\le & \frac{1}{{\nu }}\left( {\frac{{96\left| \mathcal{R} \right|{L^2}{\rho ^2}}}{\eta } + \frac{{{L^2}}}{{\left| \mathcal{R} \right|}} + \frac{{8{\rho ^2}}}{{1 - \eta }}\frac{1}{{\alpha^2}}} \right)\mathbb{E}{D_k} + \frac{{L{\sigma ^2}}}{\nu }\alpha\\
& +  \frac{4}{{\nu }} \left( {\frac{{16\left| \mathcal{R} \right| {\rho ^2}}}{\eta }\left( {{\sigma ^2} + {\zeta ^2}} \right) + n\left( {1 + \frac{{8\left| \mathcal{R} \right|{\rho ^2}}}{\eta }} \right){\varpi ^2}} \right)\\
& + \frac{1}{{\nu \alpha }}\left( {\mathbb{E}f\left( {{{\bar x}_k}} \right) - {f^*} - \left( {\mathbb{E}f\left( {{{\bar x}_{k + 1}}} \right) - {f^*}} \right)} \right).
\end{aligned}
\end{equation}
Summing (\ref{AE7-6-18}) over $k$ from 0 to $K$, $\forall K \ge 1$,  yields
\begin{equation}\label{AE7-6-19}
\begin{aligned}
&\sum\limits_{k = 0}^{K + 1} { {\mathbb{E}f\left( {{{\bar x}_k}} \right) - {f^*}} } \\
\le & \frac{4}{\nu }\left( {\frac{{16\left| \mathcal{R} \right|{\rho ^2}}}{\eta }\left( {{\sigma ^2} + {\zeta ^2}} \right) + n\left( {1 + \frac{{8\left| \mathcal{R} \right|{\rho ^2}}}{\eta }} \right){\varpi ^2}} \right)\left( {K + 1} \right) \\
& + \frac{1}{{\nu }}\left( {\frac{{96\left| \mathcal{R} \right|{L^2}{\rho ^2}}}{\eta } + \frac{{{L^2}}}{{\left| \mathcal{R} \right|}} + \frac{{8{\rho ^2}}}{{1 - \eta }}\frac{1}{{{\alpha ^2}}}} \right)\sum\limits_{k = 0}^K {\mathbb{E}{D_k}} \\
& + \frac{{\mathbb{E}f\left( {{{\bar x}_0}} \right) - {f^*}}}{{\nu \alpha }} + \frac{{L{\sigma ^2}}}{\nu }\alpha \left( {K + 1} \right).
\end{aligned}
\end{equation}
Dividing both sides of (\ref{AE7-6-19}) by $\left( {K + 1} \right)$ obtains
\begin{equation*}
\begin{aligned}
&\frac{1}{{K + 1}}\sum\limits_{k = 0}^{K + 1} {\left( {\mathbb{E}f\left( {{{\bar x}_k}} \right) - {f^*}} \right)}  \\
\le & \frac{{\mathbb{E}f\left( {{{\bar x}_0}} \right) - {f^*}}}{{\nu \alpha }\left( {K + 1} \right)} + \frac{{\frac{{96\left| \mathcal{R} \right|{L^2}{\rho ^2}}}{\eta } + \frac{{{L^2}}}{{\left| \mathcal{R} \right|}} + \frac{{8{\rho ^2}}}{{1 - \eta }}\frac{1}{{{\alpha ^2}}}}}{{\nu \alpha \left( {K + 1} \right)}}\sum\limits_{k = 0}^K {\mathbb{E}{D_k}}  \\
\end{aligned}
\end{equation*}
\begin{equation}\label{AE7-6-20}
\begin{aligned}
& + \frac{{L{\sigma ^2}}}{\nu }\alpha + {\frac{{64\left| \mathcal{R} \right|{\rho ^2}}}{{\nu \eta }}\left( {{\sigma ^2} + {\zeta ^2}} \right)} + \frac{{4n}}{\nu }\left( {1 + \frac{{8\left| \mathcal{R} \right|{\rho ^2}}}{\eta }} \right){\varpi ^2}.
\end{aligned}
\end{equation}

Recall the definition of $\mathbb{E}f_{K + 1}^{{\text{best}}}$ and then (\ref{AE7-6-20}) becomes
\begin{equation}\label{AE7-6-21}
\begin{aligned}
& \mathbb{E}f_{K + 1}^{\text{best}} - {f^*}\\
\le & \frac{{\mathbb{E}f\left( {{{\bar x}_0}} \right) - {f^*}}}{{\nu \alpha \left( {K + 1} \right)}} + \frac{{\frac{{96\left| \mathcal{R} \right|{L^2}{\rho ^2}}}{\eta } + \frac{{{L^2}}}{{\left| \mathcal{R} \right|}} + \frac{{8{\rho ^2}}}{{1 - \eta }}\frac{1}{{{\alpha ^2}}}}}{{\nu \alpha \left( {K + 1} \right)}}\sum\limits_{k = 0}^K {\mathbb{E}{D_k}}  \\
& + \frac{{L{\sigma ^2}}}{\nu }\alpha + {\frac{{64\left| \mathcal{R} \right|{\rho ^2}}}{{\eta \nu }}\left( {{\sigma ^2} + {\zeta ^2}} \right)} + \frac{{4n}}{{\nu}}\left( {1 + \frac{{8\left| \mathcal{R} \right|{\rho ^2}}}{\eta }} \right){\varpi ^2}.
\end{aligned}
\end{equation}

We then substitute (\ref{E4-1-5}) into (\ref{AE7-6-21}) and take $K$ to infinity such that (\ref{AE7-6-21}) gives rise to an asymptotic convergence error of $\mathcal{O}\left( {{\rho ^2}\left( {{\varpi ^2} + {\sigma ^2} + {\zeta ^2}} \right)} + {\varpi ^2} \right)  + \alpha \mathcal{O}\left( {{\sigma ^2}} \right) +{\alpha ^2}\mathcal{O}\left( {{\rho ^2}\left( {{\varpi ^2} + {\sigma ^2} + {\zeta ^2}} \right)} \right)$,
%\begin{equation}\label{AE7-6-22}
%\begin{aligned}
%\mathbb{E}f_{K + 1}^{\text{best}} - {f^*} \le & \mathcal{O}\left( {{\rho ^2}\left( {{\varpi ^2} + {\sigma ^2} + {\zeta ^2}} \right)} \right)  + \alpha \mathcal{O}\left( {{\sigma ^2}} \right)\\
%&+{\alpha ^2}\mathcal{O}\left( {{\rho ^2}\left( {{\varpi ^2} + {\sigma ^2} + {\zeta ^2}} \right)} \right),
%\end{aligned}
%\end{equation}
which completes the proof.
\subsection{Proof of Theorem \ref{T4}}\label{Section 7-7}
Recall the definition of $\mathcal{A}_{i,k}$ such that ${S_{i,g}} = \alpha _k \Delta_{g} $. Then, it follows from \cite[Theorem 5]{Wang2023a} that a Gaussian noise of the variance ${\varpi ^2} \ge 2\left( {\ln \left( {1.25} \right) - \ln \left( \delta  \right)} \right){\left( {{S_{i,g}}/ \epsilon} \right)^2}$ can guarantee $\left( {\epsilon, \delta } \right)$-differential privacy for any $ 0< \epsilon, \delta < 1$, which leads to (\ref{E4-1-11}) and (\ref{E4-1-12}) via substituting the upper bounds on the decaying step-size given in Theorem \ref{T2} and the constant step-size given in Theorem \ref{T3}, respectively.

\bibliographystyle{IEEEtran}
\bibliography{DP-SCC-PL}

\end{document}